# CONTOUR PROJECTED DIMENSION REDUCTION[1]

By Ronghua Luo, Hansheng Wang and Chih-Ling Tsai

*Southwestern University of Finance and Economics, Peking University and University of California, Davis*

In regression analysis, we employ contour projection (CP) to develop a new dimension reduction theory. Accordingly, we introduce the notions of the *central contour subspace* and *generalized contour subspace*. We show that both of their structural dimensions are no larger than that of the *central subspace* Cook [*Regression Graphics* (1998b) Wiley]. Furthermore, we employ CP-sliced inverse regression, CP-sliced average variance estimation and CP-directional regression to estimate the *generalized contour subspace*, and we subsequently obtain their theoretical properties. Monte Carlo studies demonstrate that the three CP-based dimension reduction methods outperform their corresponding non-CP approaches when the predictors have heavy-tailed elliptical distributions. An empirical example is also presented to illustrate the usefulness of the CP method.

**1. Introduction.** In high-dimensional data analysis, Li (1991) proposed a method of effective dimension reduction and Cook (1998b) subsequently introduced the concept of sufficient dimension reduction. Their novel approaches allow us to study low-dimensional regression relationships prior to model formulations. To effectively estimate the basis of a dimension reduction subspace, various methods have been developed. They include, but are not limited to, sliced inverse regression (SIR) [Li (1991)], sliced average variance estimation (SAVE) [Cook and Weisberg (1991)], principal Hessian directions (PHD) [Li (1992) and Cook (1998a)], minimum average variance estimator (MAVE) [Xia et al. (2002)], contour regression (CR) [Li, Zha

Received April 2008; revised December 2008.

[1]Supported in part by National Natural Science Foundation of China Grant 10771006 and a grant from Microsoft Research Asia.

and Chiaromonte (2005)], inverse regression estimation (IRE) [Cook and Ni (2005)], the Fourier method (Fourier) [Zhu and Zeng (2006)], directional regression (DR) [Li and Wang (2007)], a constructive approach [Xia (2007)], sliced regression (SR) [Wang and Xia (2008)] and also those methods based on higher-order moments [Yin and Cook (2002, 2003, 2004)].

One of the objectives of dimension reduction is to seek a *central subspace* (CS) [Cook (1994, 1998b)], which contains all information for the regression of response $Y$ on predictor $X$. To estimate the CS, two technical conditions are commonly used: the *linearity condition* [Li (1991)] and the *constant variance condition* [Cook and Weisberg (1991), Li (1992), Cook (1998a) and Li, Zha and Chiaromonte (2005)]. For example, SIR requires the *linearity condition*, while SAVE, PHD and DR entail both the *linearity* and *constant variance conditions*. It is known that the elliptically symmetric distribution of $X$ with a finite first moment implies the *linearity condition* [Li and Duan (1989)], and that the normality assumption of $X$ ensures the *constant variance condition* [Cook and Weisberg (1991)]. To facilitate the use of dimension reduction methods, Cook and Nachtsheim (1994) studied the role of elliptical symmetry in regression. In addition, they proposed a weighting procedure to achieve elliptically symmetric covariates. This motivates us to investigate dimension reduction methods via the elliptically symmetric assumption, which was also considered by Li, Zha and Chiaromonte (2005).

In the class of elliptically symmetric distributions, some either have heavy-tailed behavior or do not have finite moments. Accordingly, many existing dimension reduction methods may not yield accurate estimators of the CS. Hence, Wang, Ni and Tsai (2008) proposed the contour projection (CP) approach to project the covariate vector onto a unit contour. The resulting predictor vector has finite moments of every order and improves parameter estimators for heavy-tailed predictors. However, the theoretical properties of CP have not been thoroughly investigated. To this end, the aim of this paper is to establish a theoretical paradigm for contour projected dimension reduction. We introduce the notions of a *central contour subspace* (CCS) and a *generalized contour subspace* (GCS). Under appropriate conditions, the unique existence of the CCS and the GCS are established and their relationships with the CS are investigated. In addition, we show that their structural dimensions are no larger than that of the CS. Moreover, we obtain the theoretical properties of CP-sliced inverse regression (CP-SIR), CP-sliced average variance estimation (CP-SAVE), and CP-directional regression (CP-DR), as well as study the population exhaustiveness of those three CP methods. Consequently, the CP approach not only possesses theoretical justifications but also broadens the use of existing dimension reduction methods.

The rest of this paper is organized as follows. Section 2 introduces contour projected dimension reduction. Section 3 investigates the population features of CP-SIR and CP-SAVE, while Section 4 studies CP-DR. The sampling properties of CP methods are studied in Section 5. Extensive simulation experiments are reported in Section 6, and a real example is analyzed in Section 7. We conclude with a brief discussion in Section 8, and all technical details are left to the Appendix.

## 2. Contour projection and sufficient dimension reduction.

2.1. *Sufficient dimension reduction.* Let $X = (X_1, \ldots, X_p)^\top \in \mathbb{R}^p$ be a $p$-dimensional predictor with $p > 1$ and $Y \in \mathbb{R}^1$ be the response of interest. To capture their regression relationship, we adopt the following commonly used dimension reduction model:

$$(2.1) \qquad Y \perp\!\!\!\perp X | A^\top X,$$

where the response $Y$ is conditionally independent ($\perp\!\!\!\perp$) of $X$ given $A^\top X$ with $A \in \mathbb{R}^{p \times d}$ and $d \le p$. Let $P_A = A(A^\top A)^{-1} A^\top$ and $Q_A = I_p - P_A$, where $I_p$ is the $p$-dimensional identity matrix. As a result, the model (2.1) is equivalent to $Y \perp\!\!\!\perp X | P_A X$. For the sake of convenience, we use the generic notation $\mathcal{S}(H)$ to denote the linear subspace spanned by the column vectors of an arbitrary matrix $H$. We then refer to $\mathcal{S}(A)$ as the sufficient dimension reduction (SDR) [Cook (1998b)] subspace. When $A$ is a $p \times p$ full rank matrix, $\mathcal{S}(A)$ is automatically a SDR subspace. In practice, however, we are only interested in the "smallest" SDR subspace, which is typically defined to be the intersection of all SDR subspaces. If such an intersection is itself a SDR subspace, it is called the *central subspace* (CS) [Cook (1996, 1998b)]. Hereafter, we always assume that the CS exists and is denoted by $\mathcal{S}_{y|x}$. Next, we study the CS via CP.

2.2. *Contour projection.* For statistical validity, inverse regression methods commonly require the *linearity condition* of Li (1991), which assumes that $E(X|b^\top X)$ is a linear function of $b^\top X$, where $b \in \mathcal{S}_{y|x}$ is an arbitrary nonrandom direction. Because $b$ is unknown in practice, it is sensible to require that the *linearity condition* holds for any arbitrary direction $b \in \mathbb{R}^p$. As noted by Eaton (1986) and Cook and Nachtsheim (1994), such a requirement can only be satisfied by the so-called elliptically symmetric distribution. Its probability density function is given by [Muirhead (1982)]

$$(2.2) \qquad f_{\mu,\Sigma}(X) = |\Sigma|^{-1/2} f(\|X - \mu\|_\Sigma^2),$$

where $\mu \in \mathbb{R}^p$ is the location parameter, $\Sigma \in \mathbb{R}^{p \times p}$ is the positive definite scatter matrix and $\|X - \mu\|_\Sigma^2 = (X - \mu)^\top \Sigma^{-1} (X - \mu)$ is a Mahalanobis distance. For the sake of identifiability, we require that $\mathrm{tr}(\Sigma) = p$ [Muirhead (1982)]. Without loss of generality, we also assume that $\Sigma = I_p$ and $\mu = 0$, which can be achieved by redefining $X \doteq \Sigma^{-1/2}(X - \mu)$.

As mentioned above, if $X$ satisfies the *linearity condition* for any arbitrary direction $b$, the distribution of $X$ must be elliptically symmetric. However, this is only valid when the finite moments of an elliptically symmetric distribution exist. To avoid the issue of the existence of finite moments, we adopt the method of Wang, Ni and Tsai (2008) and propose the following contour projection operation:

$$\vec{X} = (\vec{X}_1, \ldots, \vec{X}_p)^\top = X/R, \tag{2.3}$$

where $R = \|X\|$ and $\|\cdot\|$ is the typical $L_2$ norm. As a result, $\vec{X}$ is the contour projected predictor, which has finite moments of every order. It can be shown that the support of $\vec{X}$ is the unit contour $\{\vec{x} : \|\vec{x}\| = 1\}$, as long as the support of $X$ contains an open convex set that includes the origin as an interior point. Although Wang, Ni and Tsai (2008) employed CP in the context of inverse regression, the properties of contour projected dimension reduction have not been well studied yet. This motivates us to establish the theoretical foundation for CP in the subsequent sections.

2.3. *Central contour subspace.* When $X$ follows an elliptically symmetric distribution as defined in (2.2), Wang, Ni and Tsai (2008) noted that $R$ and $\vec{X}$ are mutually independent. Accordingly, we consider the following contour projected dimension reduction model:

$$Y \perp\!\!\!\perp \vec{X} | B^\top \vec{X} \tag{2.4}$$

for some $B \in \mathbb{R}^{p \times d}$. We then label the resulting space $\mathcal{S}(B)$ a *sufficient contour subspace* (SCS) of $Y|\vec{X}$. Adopting the CS concept of Cook (1996), Cook (1998b), we define the intersection of all SCSs as the *kernel contour subspace* (KCS) and denote it by $\mathcal{K}_{y|\vec{x}}$. If $\mathcal{K}_{y|\vec{x}}$ itself is also a SCS, we call it the *central contour subspace* (CCS) and denote it by $\mathcal{C}_{y|\vec{x}}$.

Under mild yet reasonable conditions [Cook (1998b)], the CS can be well defined as the intersection of all SDR subspaces. However, in the CP context, one can easily construct examples such that the KCS is not a SCS, and hence the CCS does not exist. Consider the following example:

EXAMPLE 1.

$$Y = \sum_{j=2}^{p} X_j^2 + \varepsilon = R^2 \left( \sum_{j=2}^{p} \vec{X}_j^2 \right) + \varepsilon = R^2(1 - \vec{X}_1^2) + \varepsilon, \tag{2.5}$$

where $p > 2$ and $X$ follows an elliptically symmetric distribution. Note that $\varepsilon$ in (2.5) and hereafter satisfies $\varepsilon \perp\!\!\!\perp \vec{X}$. Let $e_j \in \mathbb{R}^p$ denote a $p$-dimensional vector with its $j$th component being 1 and others 0. Then, the second equality in (2.5) results in one SCS, $\mathcal{S}_a = \mathcal{S}(e_2, \ldots, e_p)$, while the third equality yields another SCS, $\mathcal{S}_b = \mathcal{S}(e_1)$. Nevertheless, $\mathcal{S}_a \cap \mathcal{S}_b = \varnothing$ is an empty set. Thus, the CCS does not exist. A similar example was also constructed by Cook (1994).

Example 1 indicates that the CCS may not be well defined if the regression relationship is symmetric. This motivates us to present the following three definitions so that we can assess the existence of CCS.

DEFINITION 1. We define $Y|X$ to be *dimension reducible* if (2.1) holds for some $d < p$. Otherwise, $Y|X$ is *dimension irreducible*. Analogously, we define $Y|\overrightarrow{X}$ as *dimension reducible* if (2.4) holds for some $d < p$. Otherwise, $Y|\overrightarrow{X}$ is *dimension irreducible*.

If $Y|X$ is *dimension reducible*, then (2.1) holds for $d < p$. By Lemma 2 of Wang, Ni and Tsai (2008), we have $Y \perp\!\!\!\perp \overrightarrow{X} | A^\top \overrightarrow{X}$, which implies that $Y|\overrightarrow{X}$ is also *dimension reducible*. As a result, if $Y|\overrightarrow{X}$ is *dimension irreducible*, then $Y|X$ is *dimension irreducible*. However, the reverse is not true. This indicates that the CP method might provide a better dimension reduction than SDR. To this end, we next define contour symmetric.

DEFINITION 2. Let $G_y(\overrightarrow{X} = \overrightarrow{x}) = P(Y \leq y | \overrightarrow{X} = \overrightarrow{x})$, and assume that $Y|\overrightarrow{X}$ is *dimension reducible*. We then term $Y|\overrightarrow{X}$ *contour symmetric* on directions in $\mathcal{S}(B_1)$, with $B_1 \in \mathbb{R}^{p \times d_1}$ and $1 \leq d_1 < p$, if it satisfies

$$(2.6) \quad G_y(\overrightarrow{X} = \overrightarrow{x}) = G_y(\|P_{B_1}\overrightarrow{X}\| = \|P_{B_1}\overrightarrow{x}\|, P_{B_2}\overrightarrow{X} = P_{B_2}\overrightarrow{x}),$$

where $\overrightarrow{x} \in \mathbb{R}^p$ is an arbitrary vector satisfying $\|\overrightarrow{x}\| = 1$, $B_2 \in \mathbb{R}^{p \times d_2}$ for some $0 \leq d_2 < p$ satisfying the conditions $\mathcal{S}(B_1) \cap \mathcal{S}(B_2) = \varnothing$, $\mathcal{S}(B_1) \cup \mathcal{S}(B_2) \neq \mathbb{R}^p$, and $G_y(\overrightarrow{X} = \overrightarrow{x})$ is a nondegenerate function in $\|P_{B_1}\overrightarrow{x}\|$.

Equation (2.6) implies that $\mathcal{S}(B_1, B_2)$ is a SCS. Thus, we only need to focus on the symmetric directions in SCSs. For the sake of convenience, we require that $\mathcal{S}(B_1)$ and $\mathcal{S}(B_2)$ do not have any overlap. Otherwise, one can redefine $B_1 \doteq (I_p - P_{B_2})B_1$ so that (2.6) is still valid. One might wonder why we impose the constraint $\mathcal{S}(B_1) \cup \mathcal{S}(B_2) \neq \mathbb{R}^p$. Consider an arbitrary *dimension reducible* $Y|\overrightarrow{X}$, and assume that one of its SCSs is given by $\mathcal{S}(B)$ with $\dim\{\mathcal{S}(B)\} < p$, where $\dim\{\cdot\}$ stands for the dimension of a linear subspace. If we do not impose the constraint $\mathcal{S}(B_1) \cup \mathcal{S}(B_2) \neq \mathbb{R}^p$, then we can define $B_2 = B$ and $B_1$ as a basis of linear subspace whose orthogonal complement is $\mathcal{S}(B)$. Accordingly, $\|P_{B_1}\overrightarrow{X}\|^2 = 1 - \|P_{B_2}\overrightarrow{X}\|^2$. This together with the assumption of $S(B_2)$ being a SCS implies that

$$\begin{aligned}
G_y(\overrightarrow{X} = \overrightarrow{x}) &= G_y(P_{B_2}\overrightarrow{X} = P_{B_2}\overrightarrow{x}) \\
(2.7) \quad &= G_y(\|P_{B_2}\overrightarrow{X}\| = \|P_{B_2}\overrightarrow{x}\|, P_{B_2}\overrightarrow{X} = P_{B_2}\overrightarrow{x}) \\
&= G_y(\|P_{B_1}\overrightarrow{X}\| = \|P_{B_1}\overrightarrow{x}\|, P_{B_2}\overrightarrow{X} = P_{B_2}\overrightarrow{x}).
\end{aligned}$$

Hence, condition (2.6) is satisfied for any *dimension reducible* $Y|\overrightarrow{X}$. Consequently, Definition 2 loses its ability to characterize $Y|\overrightarrow{X}$'s nontrivial symmetric structure. Therefore, requiring $\mathcal{S}(B_1) \cup \mathcal{S}(B_2) \neq \mathbb{R}^p$ is essential and necessary. In contrast to *contour symmetric*, we define *contour asymmetric* as given below.

DEFINITION 3. We call $Y|\vec{X}$ *contour asymmetric* if $Y|\vec{X}$ is not *contour symmetric* on any possible linear subspace $\mathcal{S}(B_1) \subset \mathbb{R}^p$ with $0 < \dim\{\mathcal{S}(B_1)\} < p$.

Definition 3 leads to the next theorem for the existence of the CCS.

THEOREM 1. *The central contour subspace, $\mathcal{C}_{y|\vec{x}}$, exists uniquely if and only if $Y|\vec{X}$ is contour asymmetric, that is, $\mathcal{K}_{y|\vec{x}} = \mathcal{C}_{y|\vec{x}}$.*

The above theorem provides a necessary and sufficient condition for the unique existence of the CCS, that is, the regression relationship $Y|\vec{X}$ must be *contour asymmetric*. This naturally raises an interesting question: what happens if $Y|\vec{X}$ is not *contour asymmetric*? This critical question is addressed in the next subsection.

2.4. *Generalized contour subspace.* For the existence of the CCS, Theorem 1 requires a strong condition that is likely to be violated under some symmetric situations; see Example 1 in (2.5). This motivates us to redefine the "smallest" SCS as the SCS with the smallest structural dimension, rather than the intersection of all SCSs (i.e., the KCS $\mathcal{K}_{y|\vec{x}}$). The resulting space is called the *generalized contour subspace* (GCS), denoted by $\mathcal{G}_{y|\vec{x}}$. Because $\mathbb{R}^p$ itself is a SCS and a SCS's dimension must be a positive integer no larger than $p$, there exists at least one GCS, and its dimension is unique (denoted by $d_0$). It is noteworthy that $\mathcal{K}_{y|\vec{x}} \subset \mathcal{G}_{y|\vec{x}}$, but $\mathcal{K}_{y|\vec{x}}$ is not guaranteed to be a SCS.

The existence of a GCS does not guarantee its uniqueness, so we need to find a reasonable condition to ensure the uniqueness of the GCS. To gain some insights, we consider the following example:

EXAMPLE 2.

$$(2.8) \quad Y = \vec{X}_1 + \vec{X}_2^2 + \vec{X}_3^2 + \varepsilon = \vec{X}_1 + \left(1 - \vec{X}_1^2 - \sum_{j=4}^{p} \vec{X}_j^2\right) + \varepsilon.$$

One can easily verify that there exist at least two different SCSs, whose intersection is not a SCS. For example, the first equality yields $\mathcal{G}_a = \mathcal{S}(e_1, e_2, e_3)$, while the second equality results in $\mathcal{G}_b = \mathcal{S}(e_1, e_4, \ldots, e_p)$. However, their intersection $\mathcal{G}_a \cap \mathcal{G}_b = \mathcal{S}(e_1)$ is not a SCS. If we assume that the GCS is not uniquely defined, then it is natural to have both $\mathcal{G}_a$ and $\mathcal{G}_b$ be GCSs, which implies that $\dim(\mathcal{G}_a) = \dim(\mathcal{G}_b)$. Because $\mathcal{K}_{y|\vec{x}} = \mathcal{G}_a \cap \mathcal{G}_b = \mathcal{S}(e_1)$, and $\dim(\mathcal{G}_a) = \dim(\mathcal{G}_b)$, we obtain

$$\dim(\mathcal{G}_a \setminus \mathcal{K}_{y|\vec{x}}) = \dim(\mathcal{G}_b \setminus \mathcal{K}_{y|\vec{x}}). \tag{2.9}$$

This leads to $p = 5$, so we have

$$(\mathcal{G}_a \setminus \mathcal{K}_{y|\vec{x}}) \cup (\mathcal{G}_b \setminus \mathcal{K}_{y|\vec{x}}) = (\mathbb{R}^p \setminus \mathcal{K}_{y|\vec{x}}). \tag{2.10}$$

Both (2.9) and (2.10) together imply a necessary condition for the existence of multiple GCSs. This necessary condition is

$$\begin{aligned} \dim(\mathcal{G}_a \setminus \mathcal{K}_{y|\vec{x}}) &= \tfrac{1}{2}\dim(\mathbb{R}^p \setminus \mathcal{K}_{y|\vec{x}}) \\ \iff \quad d_0 &= \tfrac{1}{2}\{p + \dim(\mathcal{K}_{y|\vec{x}})\}. \end{aligned} \tag{2.11}$$

Example 2 motivates us to find a sufficient condition for the uniqueness of the GCS. In most applications, the structural dimension $d_0$ is expected to be much smaller than the predictor dimension [Chiaromonte, Cook and Li (2002)]. Thus, a very typical violation of equality (2.11) is $d_0 < p/2 \leq \{p + \dim(\mathcal{K}_{y|\vec{x}})\}/2$. Under such a condition, the uniqueness of the GCS can be rigorously established.

THEOREM 2. *If there exists at least one SCS with structural dimension $d_0 < \{p + \dim(\mathcal{K}_{y|\vec{x}})\}/2$, then the GCS of $Y|\vec{X}$ is unique. In addition, if $\mathcal{C}_{y|\vec{x}}$ exists, then $\mathcal{G}_{y|\vec{x}} = \mathcal{C}_{y|\vec{x}}$ and $\mathcal{G}_{y|\vec{x}}$ is unique.*

The above theorem indicates that the GCS can be well defined under a rather mild condition. It also shows that the existence of the CCS implies that of the GCS. This raises another interesting question: what is the relationship between the GCS and the CS [Cook (1996) and Cook (1998b)]? This issue is addressed below.

THEOREM 3. *The relationship between the sufficient contour subspace and the CS is such that: (1) $\mathcal{K}_{y|\vec{x}} \subset \mathcal{S}_{y|x}$; and (2) $\dim(\mathcal{G}_{y|\vec{x}}) \leq \dim(\mathcal{S}_{y|x})$.*

Theorem 3 shows that the KCS is a subspace of the CS and the dimension of the GCS cannot be larger than that of the CS. In addition, Lemma 2 of Wang, Ni and Tsai (2008) indicates that the CS must be a SCS.

To further explore the relationship, one might wonder whether the GCS must be the CS. We can easily verify with the following example that this is not necessary.

EXAMPLE 3.

$$Y = \|X\|^2 + \varepsilon = R^2 + \varepsilon. \tag{2.12}$$

The first equality of Example 3 demonstrates that the dimension of the CS is the same as the predictor dimension $p$. Thus, $Y|X$ is *dimension irreducible*. However, the second equality indicates that the dimension of the GCS is only 0, which represents a substantial dimension reduction. Accordingly, the GCS is not the CS.

Since the structural dimension of the GCS is never larger than that of the CS, one might question whether the GCS is always a subspace of the CS. The answer is negative, which is illustrated by Example 1. As one can see, the first equality of (2.5) implies that the CS is $\mathcal{S}_{y|x} = \mathcal{S}(e_2, \ldots, e_p)$, while the third equality of (2.5) indicates that $\mathcal{G}_{y|\overrightarrow{x}} = \mathcal{S}(e_1)$. As a result, the intersection of the GCS and the CS is an empty set, which means GCS $\nsubseteq$ CS in this example.

In summary, Theorem 3, together with Examples 1 and 3, indicates that the GCS is closely related to but not exactly the same as the CS. Most importantly, the structural dimension of the GCS is guaranteed to be no larger, but might be much smaller, than that of the CS. Finally, one might wonder when we can have $\mathcal{G}_{y|\overrightarrow{x}} = \mathcal{S}_{y|x}$. To this end, the following theorem provides a sufficient condition.

THEOREM 4. *Assume that $Y|X$ is dimension reducible. If $Y|\overrightarrow{X}$ is contour asymmetric or $\dim(\mathcal{S}_{y|x}) < \{p + \dim(\mathcal{K}_{y|\overrightarrow{x}})\}/2$, then $\mathcal{G}_{y|\overrightarrow{x}} = \mathcal{S}_{y|x}$.*

Aforementionedly, in most applications, the structural dimension $\dim(\mathcal{S}_{y|x})$ is expected to be much smaller than the predictor dimension [Chiaromonte, Cook and Li (2002)]. Thus, the condition $\dim(\mathcal{S}_{y|x}) < \{p + \dim(\mathcal{K}_{y|\overrightarrow{x}})\}/2$ is easily satisfied as long as we have $\dim(\mathcal{S}_{y|x}) < p/2$. As a result, the technical condition entailed by Theorem 4 [i.e., $\dim(\mathcal{S}_{y|x}) < \{p + \dim(\mathcal{K}_{y|\overrightarrow{x}})\}/2$] is rather mild, which implies that the GCS is usually equal to the CS in practice. However, this condition is sufficient but not necessary. See the following counter example.

EXAMPLE 4.

$$Y = |X_1| + |X_2| + \varepsilon = R(|\overrightarrow{X}_1| + |\overrightarrow{X}_2|) + \varepsilon,$$

where $X \in \mathbb{R}^4$. In this example, we have $\mathcal{K}_{y|\overrightarrow{x}} = \varnothing$. In addition, $\dim(\mathcal{S}_{y|x}) = 2 \geq \{p + \dim(\mathcal{K}_{y|\overrightarrow{x}})\}/2 = 2$, whereas $\mathcal{S}_{y|x} = \mathcal{G}_{y|\overrightarrow{x}} = \mathcal{S}(e_1, e_2)$.

In this section, we have established the foundations of the CP, explored the properties of the GCS and the CCS, and made connections between the GCS and the CS. To facilitate the use of the GCS in dimension reduction, we now turn to studying the properties of inverse regression methods via CP in the next two sections.

**3. Contour projected SIR and SAVE.** To improve dimension reduction, Wang, Ni and Tsai (2008) employed SIR and SAVE on the CS estimation via the contour projected predictor $\overrightarrow{X}$. For the sake of simplicity, we refer to both methods as CP-SIR and CP-SAVE, respectively. Although Wang, Ni and Tsai (2008) investigated some connections between the CP approach and the CS, they failed to solve the identifiability problem due to *contour symmetry*. In addition, CP is more directly related to the GCS than the CS. These findings motivate us to establish the relationships between the linear subspaces generated by the kernel matrices of CP-SIR and CP-SAVE with the GCS, respectively.

3.1. *The CP-SIR method.* In the presence of finite moments, SIR employs the kernel matrix $\text{cov}\{E(X|Y)\}$. This naturally motivates a CP-SIR method with the kernel matrix $\text{cov}\{E(\overrightarrow{X}|Y)\}$; see Wang, Ni and Tsai (2008). Then, its statistical properties can be studied as follows.

LEMMA 1. *Assume that $X$ has the density function in (2.2). We then have*

$$\mathcal{S}\{E(\overrightarrow{X}|Y)\} \subset \mathcal{K}_{y|\overrightarrow{x}}.$$

Lemma 1 implies that CP-SIR is able to estimate a portion of the KCS $\mathcal{K}_{y|\overrightarrow{x}}$, which is a subspace of the GCS. By definition, we know that $\mathcal{K}_{y|\overrightarrow{x}} = \bigcap_{\mathcal{S}(B) \in \mathcal{I}} \mathcal{S}(B)$, where $\mathcal{S}(B)$ is a SCS and $\mathcal{I}$ is the set of all SCSs. In addition, we can decompose $\mathcal{S}(B)$ as $\mathcal{S}(B_A) \cup \mathcal{S}(B_S) = \mathcal{S}(B)$ with $\mathcal{S}(B_A) \cap \mathcal{S}(B_S) = \varnothing$ and

$$(3.1) \qquad G_y(\overrightarrow{X} = \overrightarrow{x}) = G_y(\|P_{B_S}\overrightarrow{X}\| = \|P_{B_S}\overrightarrow{x}\|, P_{B_A}\overrightarrow{X} = P_{B_A}\overrightarrow{x})$$

for some $B_A$ and $B_S$. It is noteworthy that such a decomposition always exists. This is because we can always set $S(B_s) = \varnothing$ and $\mathcal{S}(B_A)$ to be an arbitrary SCS. However, such a decomposition is not unique. Consider, for example:

EXAMPLE 5.

$$Y = |X_1| + |X_2| + X_3 + \varepsilon = R(|\overrightarrow{X}_1| + |\overrightarrow{X}_2| + \overrightarrow{X}_3) + \varepsilon,$$

where $X = (X_1, X_2, X_3) \in \mathbb{R}^3$. Here, we can decompose $\mathcal{S}(B)$ as $\{\mathcal{S}(B_S^1) = \varnothing, \mathcal{S}(B_A^1) = \mathcal{S}(e_1, e_2, e_3)\}$ (or $\{\mathcal{S}(B_S^2) = \mathcal{S}(e_1), \mathcal{S}(B_A^2) = \mathcal{S}(e_2, e_3)\}$, or $\{\mathcal{S}(B_S^3) = \mathcal{S}(e_2), \mathcal{S}(B_A^3) = \mathcal{S}(e_1, e_3)\}$). Thus, the decomposition (3.1) is not unique. In contrast, $\mathcal{K}_{y|\overrightarrow{x}} = \mathcal{S}(e_3) = \mathcal{S}(B_A^1) \cap \mathcal{S}(B_A^2) \cap \mathcal{S}(B_A^3)$ is unique, which motivates us to further characterize the properties of the KCS given below.

LEMMA 2. *Assume $X$ has the density function in (2.2). $\mathcal{S}(B)$ and $\mathcal{S}(B_A)$ are defined as above. Then: (1) $\mathcal{K}_{y|\vec{x}} = \bigcap_{\mathcal{S}(B) \in \mathcal{I}} \mathcal{S}(B_A)$; (2) If $Y|\vec{X}$ is contour symmetric on some direction, then $\mathcal{K}_{y|\vec{x}} \neq \mathcal{G}_{y|\vec{x}}$; (3) If $Y|\vec{X}$ is contour asymmetric, then $\mathcal{K}_{y|\vec{x}} = \mathcal{G}_{y|\vec{x}}$.*

Lemma 2(1) shows that the KCS excludes all symmetric directions in the SCS. Accordingly, if $Y|\vec{X}$ is *contour symmetric* on some direction, then the KCS is not a SCS, and hence cannot be a GCS. This is the result of Lemma 2(2). In addition, Lemma 2(3) indicates that the KCS becomes a GCS if $Y|\vec{X}$ is *contour asymmetric*. As a result, the KCS becomes the CCS (see Theorem 1). This finding bridges $\mathcal{S}\{E(\vec{X}|Y)\}$ and the KCS by introducing the following assumption.

ASSUMPTION 1. *For any $v \in \mathcal{K}_{y|\vec{x}}$, $v \neq 0$, $E(v^\top \vec{X}|Y)$ is nondegenerate.*

Assumption 1 is similar to the mild requirement given in Li and Wang (2007) and Shao, Cook and Weisberg (2007). Under this assumption and the *contour asymmetric* condition, we establish the population exhaustiveness of CP-SIR, given below.

THEOREM 5. *Assume $X$ has the density function in (2.2). If Assumption 1 holds and $Y|\vec{X}$ is contour asymmetric, then $\mathcal{S}\{E(\vec{X}|Y)\} = \mathcal{G}_{y|\vec{x}}$.*

To facilitate the use of CP, we next study CP-SAVE.

3.2. *The CP-SAVE method.* As demonstrated in the previous subsection, CP-SIR fails if the regression relationship is symmetric. Under such a situation, SAVE [Cook and Weisberg (1991)] may provide a better approach to dimension reduction. Thus, it is of great interest to explore the usefulness of SAVE with contour projected predictors.

For the sake of simplicity, we denote $\mathcal{G}_{y|\vec{x}} = \mathcal{S}(B_0)$, where $B_0 \in \mathbb{R}^{p \times d_0}$ is an orthonormal basis. Let $\tau(Y) = \{1 - \lambda(Y)\}/(p - d_0)$ and $\lambda(Y) = E(\|B_0^\top \vec{X}\|^2 | Y)$. Then, it can be verified that $Q_{B_0} E(\vec{X} \vec{X}^\top | Y) Q_{B_0} = \tau(Y) Q_{B_0}$, and the estimator of $\tau(Y)$ is given in Section 5.4. This motivates us to consider the CP-SAVE kernel matrix, $M_{\text{SAVE}} = E\{[\tau(Y) I_p - E(\vec{X} \vec{X}^\top | Y)]^2\}$, where $E(\vec{X} \vec{X}^\top | Y) = \text{cov}(\vec{X}|Y)$. The following lemma shows that CP-SAVE enables us to estimate a portion of the GCS.

LEMMA 3. *Under (2.2), $\mathcal{S}\{\tau(Y) I_p - E(\vec{X} \vec{X}^\top | Y)\} \subset \mathcal{G}_{y|\vec{x}}$.*

It is known that the traditional SAVE method is able to estimate a portion of the CS under both the *linearity condition* and the *constant variance* *condition*. Lemma 3 demonstrates that CP-SAVE can estimate a portion of the GCS when the predictor is elliptically symmetric distributed. It is remarkable that the *constant variance condition* is no longer needed here. Furthermore, as noted by Cook and Lee (1999) and Li and Wang (2007), SAVE estimates the CS exhaustively under some reasonable conditions. This motivates us to study whether CP-SAVE can estimate the GCS exhaustively. To this end, we need the following assumption.

ASSUMPTION 2. Let $\beta_{0i}$ be the $i$th component of $B_0$ and $w = (w_1, \ldots, w_{d_0})^\top$ be an arbitrary $d_0 \times 1$ nonzero vector, where $B_0$ is defined as above. We then assume that $\sum_{i=1}^{d_0} w_i \phi(Y, \beta_{0i})$ is a nondegenerate random variable with $\phi(Y, \beta_{0i}) = E\{(\beta_{0i}^\top \overrightarrow{X})^2 | Y\}$.

Assumption 2 is valid only if $Y|\overrightarrow{X}$ is *dimension reducible*. Otherwise, we have $\sum_{i=1}^{p} \phi(Y, \beta_{0i}) = 1$, which clearly violates Assumption 2. A similar assumption on the predictor $X$ can be found in Li and Wang (2007) and Shao, Cook and Weisberg (2007). With the help of Assumption 2, the population exhaustiveness of CP-SAVE can be established.

THEOREM 6. *Assume $X$ has the density function in (2.2), and $Y|\overrightarrow{X}$ is dimension reducible. Then, under Assumption 2, we have*

$$\mathcal{S}\{\tau(Y)I_p - E(\overrightarrow{X}\overrightarrow{X}^\top | Y)\} = \mathcal{G}_{y|\overrightarrow{x}}.$$

In addition to CP-SIR and CP-SAVE, we further extend the CP approach to the dimension reduction method proposed by Li and Wang (2007) in the next section.

## 4. Contour projected DR.

4.1. *Motivation*. Both SIR and SAVE are commonly used dimension reduction estimators. However, SIR fails to provide exhaustive estimation for symmetric regressions [Li (1991) and Cook and Weisberg (1991)]. Although SAVE is able to estimate the CS exhaustively [Cook and Lee (1999) and Li and Wang (2007)], its estimation efficiency is relatively poor. To mitigate the weaknesses and enhance the strengths of SIR and SAVE, Li and Wang (2007) recently proposed a novel method, directional regression (DR), for dimension reduction. Specifically, they suggested the kernel matrix of DR to be $E\{2I_p - A^0(Y, Y^*)\}^2$, where $A^0(Y, Y^*) = E\{(X - X^*)(X - X^*)^\top | Y, Y^*\}$, and $(X^*, Y^*)$ is an independent copy of $(X, Y)$. Moreover, Li and Wang (2007) demonstrated that this matrix is closely related to that of SIR and SAVE in a very interesting yet effective manner.

The major advantages of DR are: (i) DR achieves higher estimation efficiency, and requires fewer computations; (ii) DR is able to estimate the CS exhaustively. These nice properties motivate us to propose a dimension reduction method of contour projected directional regression (CP-DR) that synthesizes the strengths from both the CP and DR approaches. To this end, we adopt Li and Wang's (2007) approach to consider the following kernel matrix for CP-DR:

$$M_{\mathrm{DR}} = E[\{\tau(Y) + \tau(Y^*)\}I_p - A(Y, Y^*)]^2,$$

where $A(Y, Y^*) = E\{(\overrightarrow{X} - \overrightarrow{X}^*)(\overrightarrow{X} - \overrightarrow{X}^*)^\top | Y, Y^*\}$, and $(\overrightarrow{X}^*, Y^*)$ is an independent copy of $(\overrightarrow{X}, Y)$. Analogously to Li and Wang (2007), we refer to $\overrightarrow{X} - \overrightarrow{X}^*$ as the contour projected empirical direction. Note that, due to the absence of the *constant variance condition*, we replace the constant 2 in the DR kernel matrix with $\{\tau(Y) + \tau(Y^*)\}$ to constitute the CP-DR kernel matrix.

The CP-DR dimension reduction method inherits all nice properties from DR. Furthermore, CP-DR enables us to handle heavy-tailed predictor distributions. Finally, CP-DR has the potential to produce a dimension reduction subspace with a much smaller structural dimension than that of DR (see Example 3).

4.2. *Population exhaustiveness.* Similar to CP-SIR and CP-SAVE, we present the following result to assure that CP-DR estimates a portion of the GCS.

THEOREM 7. *Under (2.2), $\mathcal{S}\{[\tau(Y) + \tau(Y^*)]I_p - A(Y, Y^*)\} \subset \mathcal{G}_{y|\overrightarrow{x}}$.*

Because DR estimates the CS exhaustively, it is of great interest to show that CP-DR can also estimate the GCS exhaustively.

THEOREM 8. *Assume that either: (1) Assumption 1 holds and $Y|\overrightarrow{X}$ is contour asymmetric; or (2) Assumption 2 holds and $Y|\overrightarrow{X}$ is dimension reducible. Under (2.2), we then have $\mathcal{S}(M_{\mathrm{DR}}) = \mathcal{G}_{y|\overrightarrow{x}}$.*

Theorem 8 indicates that CP-DR estimates the GCS exhaustively without employing the *constant variance condition*, which is used by various exhaustive methods [Li, Zha and Chiaromonte (2005), Zhu and Zeng (2006) and Li and Wang (2007)].

4.3. *A simplified formulation.* To reduce the computations necessary for CP-DR, we simplify the analytical form of the kernel matrix $M_{\mathrm{DR}}$ in the next theorem.

THEOREM 9. *The matrix $M_{\mathrm{DR}}$ can be expressed as*

$$M_{\mathrm{DR}} = 2[E\{\tau^2(Y)\}I_p + E\{E^2(\vec{X}\vec{X}^\top|Y)\} + E^2\{E(\vec{X}|Y)E(\vec{X}^\top|Y)\}$$
(4.1)
$$+ E\{E(\vec{X}^\top|Y)E(\vec{X}|Y)\}E\{E(\vec{X}|Y)E(\vec{X}^\top|Y)\}$$
$$- 2E\{\tau(Y)E(\vec{X}\vec{X}^\top|Y)\}].$$

To ease interpretation, one can rewrite (4.1) as

$$M_{\mathrm{DR}} = 2[E\{E^2[\tau(Y)I_p - \vec{X}\vec{X}^\top|Y]\} + E^2\{E(\vec{X}|Y)E(\vec{X}^\top|Y)\}$$
$$+ E\{E(\vec{X}^\top|Y)E(\vec{X}|Y)\}E\{E(\vec{X}|Y)E(\vec{X}^\top|Y)\}].$$

Thus, it is a natural combination of the kernel matrices of CP-SAVE and CP-SIR. According to Theorem 9, we are able to estimate $\mathcal{G}_{y|\vec{x}}$ as long as we can estimate $\tau(Y)$, $E(\vec{X}\vec{X}^\top|Y)$ and $E(\vec{X}|Y)$ consistently. The parameter estimators and their properties are discussed in the next section.

## 5. The sampling properties.

5.1. *The estimators of $\mu$, $\Sigma$ and $\vec{x}$.* Without loss of generality, we theoretically assume that $\mu = 0$ and $\Sigma = I_p$ (see Section 2.2). In practice, however, both $\mu$ and $\Sigma$ are often unknown and have to be estimated from the data. It is noteworthy that both the distributions of $X$ and $\vec{X}$ are elliptically symmetric with the same contour shapes. Consequently, they share the same scatter matrix $\Sigma$ via the identifiable constraint $\mathrm{tr}(\Sigma) = p$. Thus, the estimated scatter matrix of $\vec{X}$ can be used to estimate that of $X$; see Tyler (1987).

More specifically, let $(y_i, x_i)$ be the observation collected from the $i$th subject $(1 \leq i \leq n)$, where $y_i \in \mathbb{R}^1$ is the response and $x_i = (x_{i1}, \ldots, x_{ip})^\top \in \mathbb{R}^p$ is a $p \times 1$ predictor vector. To estimate $\mu$ and $\Sigma$, we follow the method of Tyler (1987) and define $\hat{\mu}^{(0)} = (\hat{\mu}_1^{(0)}, \ldots, \hat{\mu}_p^{(0)})^\top$, where $\hat{\mu}_j^{(0)}$ is the median of $\{x_{ij} : i = 1, \ldots, n\}$ for every $1 \leq j \leq p$. As a result, $\hat{\mu}^{(0)}$ is $\sqrt{n}$-consistent. Next, we estimate $\Sigma$ and $\mu$ by iterating the following two equations [Tyler (1987)]:

$$(5.1) \qquad \hat{\Sigma}^{(m+1)} \propto n^{-1} \sum_{i=1}^n \frac{(x_i - \hat{\mu}^{(m)})(x_i - \hat{\mu}^{(m)})^\top}{\|x_i - \hat{\mu}^{(m)}\|^2_{\hat{\Sigma}^{(m)}}}$$

and

$$(5.2) \quad \hat{\mu}^{(m+1)} = \left(n^{-1}\sum_{i=1}^n \frac{1}{\|x_i - \hat{\mu}^{(m)}\|_{\hat{\Sigma}^{(m+1)}}}\right)^{-1} \left(n^{-1}\sum_{i=1}^n \frac{x_i}{\|x_i - \hat{\mu}^{(m)}\|_{\hat{\Sigma}^{(m+1)}}}\right),$$

where $\|\cdot\|_\Sigma$ is the Mahalanobis norm as defined in Section 2.2 and $(\hat{\mu}^{(m)}, \hat{\Sigma}^{(m)})$ is the estimator of $(\mu, \Sigma)$ obtained in the $m$th step. We iterate the procedures of (5.1) and (5.2) until they converge, and denote the resulting estimator by $(\hat{\mu}, \hat{\Sigma})$. For the sake of simplicity, one can fix $\hat{\mu}^{(m)} = \hat{\mu}^{(0)}$ (for every $m \geq 1$) without iterating (5.2). The asymptotic efficiency of $\hat{\Sigma}$ is not affected as long as the predictor distribution is elliptically symmetric. Then, by Tyler's (1987) Theorem 2.2, this iterating process is guaranteed to converge computationally with probability 1. Furthermore, by Tyler's (1987) Theorem 4.2, the estimator $\hat{\Sigma}$ is $\sqrt{n}$-consistent. In other words, we have $\|\hat{\Sigma} - \Sigma\| = O_p(n^{-1/2})$, where $\|H\|$ is defined to be the maximum of the absolute singular value of a matrix $H$, and it becomes the usual $L_2$ norm if $H$ is a vector. Using $\hat{\mu}$ and $\hat{\Sigma}$, we obtain the estimator of the contour projected predictor, $\overrightarrow{x}_i = \hat{\Sigma}^{-1/2}(x_i - \hat{\mu})/\|x_i - \hat{\mu}\|_{\hat{\Sigma}}$.

5.2. *A preliminary result.* To establish the $\sqrt{n}$-consistency of the three CP estimators, we first introduce a technical assumption and then present one preliminary result.

ASSUMPTION 3. Assume $\|X\|^2$ has a continuous distribution with a probability density function $h(\cdot)$. We further assume that there exist constants $\alpha > 1$ and $C_\alpha > 0$ such that $t^{-\alpha} h(t) \to C_\alpha$ as $t \to 0$.

This is a very mild yet reasonable assumption. For example, if $X$ follows a standard normal distribution, then $\|X\|^2$ is a chi-square distribution with $p$ degrees of freedom. Thus, Assumption 3 is satisfied with $\alpha = p/2 - 1$ as long as $p > 4$. Analogously, if $X$ follows a multivariate $t$-distribution with $df$ degrees of freedom, then $\|X\|^2$ is a $F$-distribution with the degrees of freedom $(p, df)$. Once again, Assumption 3 is satisfied with $\alpha = p/2 - 1$ as long as $p > 4$. Applying the above assumption, we obtain the following preliminary result.

THEOREM 10. *Under Assumption 3, we have:* (i) $E\|X\|^{-4} < \infty$. *Furthermore, assume that $\tilde{\mu} \in \mathbb{R}^p$ is an arbitrary random vector satisfying $\tilde{\mu} = O_p(n^{-1/2})$, and $\tilde{\Sigma} \in \mathbb{R}^{p \times p}$ is an arbitrary random matrix satisfying $\tilde{\Sigma}_p \to I_p$. We then have:*

(ii) $$\max_{1 \leq i \leq n} \left| \frac{\|x_i - \tilde{\mu}\|_{\tilde{\Sigma}}^2}{\|x_i\|_{\tilde{\Sigma}}^2} - 1 \right| \to_p 0.$$

Theorem 10 (i) indicates that $\|X\|^{-4}$ has a finite moment, which plays an important role in ensuring the $\sqrt{n}$-consistency of $\hat{\Sigma}$; see Theorem 4.2 of Tyler (1987). Hence, it will be useful for us to show the $\sqrt{n}$-consistency of the CP estimators. In addition, Theorem 10(ii) allows us to replace $\|x_i - \tilde{\mu}\|_{\tilde{\Sigma}}$ by $\|x_i\|_{\tilde{\Sigma}}$, which simplifies the theoretical proof of $\sqrt{n}$-consistency; see Appendix A.14 for the details.

5.3. *$\sqrt{n}$-consistency of CP-SIR.* Without loss of generality, we assume that $Y$ is discrete and has a finite support $\{1,\ldots,K\}$ with $K \geq 2$ [see Li (1991) and Cook and Ni (2005)]. Next, we define $z_{ik} = 1$ if $Y = k$ and otherwise 0. Then, $n_k = \sum_i z_{ik}$ is the number of observations falling into the $k$th slice. Under this setting, the kernel matrix of CP-SIR can be defined as $M_{\text{SIR}} = \text{cov}\{E(\overrightarrow{X}|Y)\} = \Gamma_{\text{SIR}}\Gamma_{\text{SIR}}^\top$, where $\Gamma_{\text{SIR}} = \{E(\overrightarrow{X}|Y=1)\sqrt{p_1},\ldots, E(\overrightarrow{X}|Y=K)\sqrt{p_K}\} \in \mathbb{R}^{p \times K}$ and $p_k = P(Y=k) = E(z_{ik})$. Thus, it is natural to estimate $M_{\text{SIR}}$ by $\hat{M}_{\text{SIR}} = \hat{\Gamma}_{\text{SIR}}\hat{\Gamma}_{\text{SIR}}^\top$, where $\hat{\Gamma}_{\text{SIR}} = (\bar{x}_1\sqrt{\hat{p}_1},\ldots,\bar{x}_K\sqrt{\hat{p}_K})$, $\hat{p}_k = n_k/n$, and $\bar{x}_k = n_k^{-1}\sum_i \overrightarrow{x}_i z_{ik}$. In the following theorem, we show that $\hat{M}_{\text{SIR}}$ is $\sqrt{n}$-consistent.

THEOREM 11. *Under Assumption 3, we have $\bar{x}_k - E\{\overrightarrow{X}|Y=k\} = O_p(n^{-1/2})$.*

The above theorem together with $\hat{p}_k - p_k = O_p(n^{-1/2})$ implies $\|\hat{M}_{\text{SIR}} - M_{\text{SIR}}\| = O_p(n^{-1/2})$. Hence, CP-SIR achieves $\sqrt{n}$-consistency.

5.4. *$\sqrt{n}$-consistency of CP-SAVE and CP-DR.* Under the assumption that $Y$ is discrete, the kernel matrix of CP-SAVE given by $M_{\text{SAVE}} = E[\tau(Y)I_p - E(\overrightarrow{X} \times \overrightarrow{X}^\top|Y)]^2$ in Section 3.2 can be written as $M_{\text{SAVE}} = \sum_k p_k E[\tau_k I_p - E(\overrightarrow{X}\overrightarrow{X}^\top|Y=k)]^2$, where $\tau_k = \tau(Y=k)$. For the given data, we can estimate this kernel matrix by $\hat{M}_{\text{SAVE}} = \sum_k \hat{p}_k(\hat{\tau}_k I_p - \hat{\Sigma}_k)^2$, where $\hat{\Sigma}_k = n_k^{-1}\sum_i \overrightarrow{x}_i \overrightarrow{x}_i^\top z_{ik}$ and $\hat{\tau}_k$ is the median of $\hat{\Sigma}_k$'s eigenvalues. The reason for using the median of $\hat{\Sigma}_k$'s eigenvalues to estimate $\tau_k$ is as follows. By Theorem 6, we know that most of the eigenvalues of $\Sigma_k = E(\overrightarrow{X}\overrightarrow{X}^\top|Y=k)$ are equal to $\tau_k$, except for those eigenvalues associated with $\mathcal{G}_{y|\overrightarrow{x}}$. As a result, the median of $\Sigma_k$'s eigenvalues is $\tau(Y)$ whenever $d_0 < p/2$. Applying similar techniques to those used in the proof of Theorem 11, we are able to show that $\|\hat{\Sigma}_k - \Sigma_k\| = O_p(n^{-1/2})$. Consequently, we have $\hat{\tau}_k - \tau_k = O_p(n^{-1/2})$; see Eaton and Tyler (1994). Therefore, $\|\hat{M}_{\text{SAVE}} - M_{\text{SAVE}}\| = O_p(n^{-1/2})$, which implies that CP-SAVE is $\sqrt{n}$-consistent. Furthermore, by Theorem 9, we find that $M_{\text{DR}}$ is closely related to those of CP-SIR and CP-SAVE. Therefore, $M_{\text{DR}}$ can be estimated analogously, and the resulting estimator is also $\sqrt{n}$-consistent.

5.5. *The estimation of structural dimension.* In this subsection, we propose an informal but effective method for determining the structural dimension. For the sake of convenience, we use the generic notation $M$ and $\hat{M}$ to represent a kernel matrix and its consistent estimator, respectively. In addition, we assume that the structure dimension of $M$ is $d_0$, that is, $\lambda_j > 0$ for any $j \leq d_0$ but $\lambda_j = 0$ for $j > d_0$, where $\lambda_j$ is the $j$th largest eigenvalue of $M$. Moreover, we define $\hat{r}_j = \hat{\lambda}_j/\hat{\lambda}_{j+1}$ for $1 \leq j \leq p-1$. Intuitively, if $j < d_0$, both $\hat{\lambda}_j$ and $\hat{\lambda}_{j+1}$ converge in probability to positive constants. Thus, we have $\hat{r}_j = O_p(1)$ for any $j < d_0$. On the other hand, if $j > d_0$, both $\hat{\lambda}_j$ and $\hat{\lambda}_{j+1}$ converge in probability to 0. We then assume that both $\hat{\lambda}_j$ and $\hat{\lambda}_{j+1}$ share the same convergence speed [this assumption is indeed valid if $\sqrt{n}(\hat{M} - M)$ is asymptotically normal; see Eaton and Tyler (1994)]. Under this assumption, we also have $\hat{r}_j = O_p(1)$ for any $j > d_0$. However, if $j = d_0$, then $\hat{r}_j \to \infty$. This is because $\hat{\lambda}_{d_0} \to \lambda_{d_0} > 0$ but $\hat{\lambda}_{d_0+1} \to 0$. Consequently, one can estimate $d_0$ by $\hat{d} = \arg\max_{1 \le j \le d_{\max}} \hat{r}_j$, where $d_{\max}$ is the maximum dimension given a priori. In practice, we recommend using $d_{\max} = 5$, which is large enough for the purpose of dimension reduction. We refer to this estimation method as Maximal Eigenvalue Ratio Criterion (MERC), and simulation studies in Section 6.3 suggest that such a simple method works fairly well.

## 6. Monte Carlo studies.

6.1. *Simulation settings.* To evaluate the finite sample performance of CP methods, we conducted extensive Monte Carlo simulations. In these studies, $\tilde{X} = (X^\top, \varepsilon)^\top$ were independently generated from $W/\sqrt{V_{df}/df}$, where $W \in \mathbb{R}^{p+1}$ is a $(p+1)$-dimensional standard normal random variable, $V_{df}$ is a chi-squared random variable with degrees of freedom $df$, and $V_{df}$ is independent of $W$. As a result, $X$ follows a multivariate $t$ distribution with $df$ degrees of freedom [see Lange, Little and Taylor (1989)]. In addition, the marginal distributions of $X$ and $\varepsilon$ are $p$-multivariate $t$ distribution and univariate $t$ distribution, respectively [see Muirhead's (1982)]. Thus, the probability density function of $X$ has the form of (2.2). We then simulated $Y = g(B_0^\top \overrightarrow{X}, \|X\|, \varepsilon)$, where $B_0 = (\beta_{01}^\top, \ldots, \beta_{0d_0}^\top)^\top$ and $g(\cdot, \cdot, \cdot)$ are the pre-specified functions of each model given in the next subsection. Under this setting, one can show that $\overrightarrow{X}$ is independent of $(\|X\|, \varepsilon)$, which leads to $Y \perp\!\!\!\perp \overrightarrow{X} | B_0^\top \overrightarrow{X}$. It is noteworthy that $X$ and $\varepsilon$ are not independent, except when $X$ follows a normal distribution (i.e., $df = \infty$).

There are five models considered in our simulation studies. For each model, we simulated 500 data sets. In addition, the number of covariates is $p = 20$ and the number of slices is $H = 5$. Moreover, four degrees of freedoms are examined, namely $df = 1, 3, 5$ and $\infty$. They represent the case that moments do not exist, the first moment exists, the second moment exists, and the distribution is normal, respectively. For each model, three different inverse regression methods (i.e., SIR, SAVE and DR) and their corresponding contour projected approaches are compared. To evaluate the accuracy of the estimate $\hat{B}$, we adopt Li and Wang's (2007) distance measure, $\Delta(B_0, \hat{B}) = \text{tr}\{(P_{B_0} - P_{\hat{B}})(P_{B_0} - P_{\hat{B}})^\top\}/d_0$. A smaller value of $\Delta$ indicates a better estimate. For the sake of simplicity, we slightly abuse notation by using $B_0 = (\beta_{01}^\top, \ldots, \beta_{0d_0}^\top)^\top$ to commonly represent the true parameter matrix for each of the five models.

6.2. *Estimation accuracy.* To evaluate the estimation accuracy of the proposed CP methods, we present the detailed structures of the five models given below:

(I) LINEAR CONDITIONAL MEAN MODEL [Li (1991) and Ni, Cook and Tsai (2005)], $Y = \beta_{01}^\top X + 0.5\varepsilon = \beta_{01}^\top \overrightarrow{X} \times \|X\| + 0.5\varepsilon$, where $\beta_{01} = (1,1,1,0,\ldots,0)^\top \in \mathbb{R}^p$ and $d_0 = 1$.

(II) SYMMETRIC CONDITIONAL MEAN MODEL [Li, Zha and Chiaromonte (2005)], $Y = (\beta_{01}^\top X_1)^2 + \beta_{02}^\top X_2 + 0.2\varepsilon = (\beta_{01}^\top \overrightarrow{X}_1)^2 \times \|X\|^2 + \beta_{02}^\top \overrightarrow{X}_2 \times \|X\| + 0.2\varepsilon$, where $\beta_{01} = (1,0,\ldots,0)^\top \in \mathbb{R}^p$, $\beta_{02} = (0,1,\ldots,0)^\top \in \mathbb{R}^p$ and $d_0 = 2$.

(III) DISCRETE RESPONSE MODEL [Zhu and Zeng (2006)], $Y = I(\beta_{01}^\top X + 0.2\varepsilon > 0) + 2I(\beta_{02}^\top X + 0.2\varepsilon > 0) = I(\beta_{01}^\top \overrightarrow{X} \times \|X\| + 0.2\varepsilon > 0) + 2I(\beta_{02}^\top \overrightarrow{X} \times \|X\| + 0.2\varepsilon > 0)$, where $I(\cdot)$ denotes the indicator function, $\beta_{0j} \in \mathbb{R}^p$ ($j=1,2$), and $d_0 = 2$. The first four components of $\beta_{01}$ and the seventh to tenth components of $\beta_{02}$ are taken to be 1, while the rest of the components of $\beta_{01}$ and $\beta_{02}$ are fixed to be 0.

(IV) HETEROGENEOUS VARIANCE MODEL [Li, Zha and Chiaromonte (2005)], $Y = 0.5(\beta_{01}^\top X - 0.5)^2 \varepsilon = 0.5(\beta_{01}^\top \overrightarrow{X}\|X\| - 0.5)^2 \varepsilon$, where $\beta_{01} = (1,0,\ldots,0)^\top \in \mathbb{R}^p$ and $d_0 = 1$.

(V) CONTOUR SYMMETRIC MODEL, $Y = (\beta_{01}^\top X)^2 \times \|X\|^{-2} + 0.2\varepsilon = (\beta_{01}^\top \times \overrightarrow{X})^2 + 0.2\varepsilon$, where $\beta_{01} = (1,1,0,\ldots,0)^\top \in \mathbb{R}^p$ and $d_0 = 1$.

It is noteworthy that GCS = CS in models I to IV, while GCS $\neq$ CS in model V. We consider three sample sizes ($n = 200$, 400 and 1000) in Monte Carlo studies. To save space, Table 1 only reports the results with $n = 400$. We find that the performances of the three non-CP methods deteriorates seriously as the tail of the predictor distribution becomes heavier. However, the CP-methods (particularly the CP-DR method) perform satisfactorily across all simulation settings.

6.3. *Dimension determination.* Employing models I–V discussed in the previous subsection, we study the finite sample performance of MERC. Because the three CP methods and the four *df*'s yield qualitatively similar findings, we only report the results of CP-DR with $df = 3$. Table 2 indicates that the percentage of $\hat{d} = d_0$ steadily approaches 100% as the sample size increases. Consequently, MERC is a simple and effective method to determine the structural dimension in large samples.

TABLE 1
*The average of $\Delta(B_0, \hat{B})$ for models* I–V *with* $n = 400$

| Model | $df$ | CP-DR | DR | CP-SIR | SIR | CP-SAVE | SAVE |
|---|---|---|---|---|---|---|---|
| I | $\infty$ | 0.024 | 0.023 | 0.020 | 0.020 | 0.096 | 0.096 |
|   | 5 | 0.030 | 0.062 | 0.026 | 0.058 | 0.130 | 0.810 |
|   | 3 | 0.033 | 0.517 | 0.028 | 0.137 | 0.158 | 1.693 |
|   | 1 | 0.580 | 1.859 | 0.076 | 0.800 | 1.238 | 1.817 |
| II | $\infty$ | 0.126 | 0.130 | 0.989 | 0.994 | 0.313 | 0.314 |
|   | 5 | 0.160 | 0.678 | 0.992 | 1.030 | 0.448 | 0.965 |
|   | 3 | 0.183 | 1.132 | 0.993 | 1.084 | 0.537 | 1.172 |
|   | 1 | 0.828 | 1.526 | 1.035 | 1.492 | 1.130 | 1.528 |
| III | $\infty$ | 0.169 | 0.177 | 0.192 | 0.194 | 0.650 | 0.699 |
|   | 5 | 0.165 | 0.678 | 0.185 | 0.233 | 0.643 | 1.427 |
|   | 3 | 0.172 | 1.294 | 0.190 | 0.316 | 0.639 | 1.674 |
|   | 1 | 0.591 | 1.772 | 0.214 | 0.973 | 1.087 | 1.789 |
| IV | $\infty$ | 0.238 | 0.244 | 0.473 | 0.513 | 0.397 | 0.370 |
|   | 5 | 0.317 | 1.198 | 0.600 | 0.925 | 0.495 | 1.256 |
|   | 3 | 0.409 | 1.719 | 0.736 | 1.334 | 0.573 | 1.726 |
|   | 1 | 1.534 | 1.878 | 1.252 | 1.880 | 1.587 | 1.879 |
| V | $\infty$ | 0.302 | 0.341 | 1.897 | 1.900 | 0.302 | 0.341 |
|   | 5 | 0.401 | 1.575 | 1.894 | 1.894 | 0.401 | 1.575 |
|   | 3 | 0.464 | 1.839 | 1.891 | 1.886 | 0.463 | 1.839 |
|   | 1 | 1.549 | 1.896 | 1.889 | 1.906 | 1.548 | 1.896 |

6.4. *Asymmetric predictors.* In this paper, the CP theory requires the distribution of the predictor vector to be elliptically symmetric. However, in practice, such an assumption could be violated to some extent. Thus, it is of interest to evaluate the CP performance under nonelliptically symmetric distributions. To this end, we regenerate the $W$ random vector (defined in Section 6.1) from a centralized standard exponential distribution, that is, $\exp(1) - 1$. We then replicate the same simulation experiments as given in the previous subsection. Because the results are qualitatively similar to those in Table 1, we only report the results with $df = 3$ and $n = 400$. Table 3

TABLE 2
*The percentage $\hat{d} = d_0$ for CP-DR with* $df = 3$

|  | The five models | | | | |
|---|---|---|---|---|---|
| $n$ | I | II | III | IV | V |
| 200 | 0.996 | 0.588 | 0.324 | 0.620 | 0.658 |
| 400 | 1.000 | 0.962 | 0.642 | 0.916 | 0.904 |
| 1000 | 1.000 | 1.000 | 1.000 | 1.000 | 1.000 |

TABLE 3
*The average of $\Delta(B_0, \hat{B})$ under asymmetric predictors with $(n, df) = (400, 3)$*

| Model | CP-DR | DR | CP-SIR | SIR | CP-SAVE | SAVE |
|---|---|---|---|---|---|---|
| I | 0.049 | 0.682 | 0.034 | 0.154 | 0.737 | 1.764 |
| II | 0.092 | 0.691 | 0.542 | 0.827 | 0.164 | 0.822 |
| III | 0.153 | 1.317 | 0.130 | 0.254 | 1.164 | 1.717 |
| IV | 0.719 | 1.675 | 1.307 | 1.813 | 0.781 | 1.677 |
| V | 0.289 | 1.859 | 1.343 | 1.464 | 0.298 | 1.862 |

shows that the CP methods perform reasonably well even with asymmetric predictors.

**7. A real example.** To demonstrate the practical usefulness of the CP methodology, we consider an example from the Chinese stock market. The dataset is obtained from the CCER database, which is one of the most authoritative commercial databases for the Chinese stock market (<http://www.ccerdata.com/>). It contains yearly accounting information for the firms that are publicly listed in the Chinese stock market during the period from 1997 to 2000. The total sample size is 2951 observations. The objective of this study is to understand these firms' earnings patterns, which can be useful information for an investment decision. To this end, the response is the firm's next year return on equity (ROEt). The predictors include the current year accounting variables: return on equity (ROE), log-transformed total assets (ASSET), profit margin ratio (PM), sales growth rate (GROWTH), leverage level (LEV) and asset turnover ratio (ATO). A simple calculation shows that the averaged yearly kurtosis of the aforementioned explanatory variables are given by 260.66 (ROE), 2.82 (ASSET), 365.03 (PM), 224.98 (GROWTH), 10.04 (LEV) and 10.64 (ATO), respectively. Thus, all predictors other than ASSET have heavy-tailed distributions, which motivates us to employ CP for estimating parameters. All predictors are appropriately scaled so that the resulting diagonal components of the scatter matrix are equal to one.

Because CP-DR yields more reliable estimates than the other methods in simulation studies, we employ CP-DR to analyze the data. The resulting structural dimension estimated by MERC is $\hat{d} = 2$, and the first two CP-DR estimates are given below.

| Direction | ROE | ASSET | PM | GROWTH | LEV | ATO |
|---|---|---|---|---|---|---|
| $\hat{\beta}_1$ | 0.936 | −0.008 | 0.120 | 0.139 | 0.113 | 0.280 |
| $\hat{\beta}_2$ | −0.279 | −0.045 | 0.825 | 0.167 | 0.168 | 0.428 |

The above estimates clearly indicate that ROE, PM and ATO are the most important variables associated with the firm's future earnings. To further understand their effects, Figure 1 depicts the scatter plots of $Y$ versus $\hat{\eta}_j = \hat{\beta}_j^\top \overrightarrow{X}$ $(j = 1, 2)$. For a better view, we trim off an extremity of 5% of the observations according to the response ROEt. As a result, both plots display the monotonically increasing pattern for those observations with ROEt $> 0$. Practically, the first index $\hat{\eta}_1$ can be viewed as the autocorrelation index because the majority of the weight is loaded on the predictor ROE. Moreover, the second index $\hat{\eta}_2$ can be regarded as the ratio index because the notable weights are loaded on the ratio predictors PM and ATO. Consequently, the higher the return on equity or the larger the profit margin ratio and turnover ratio, the greater the yield of future return (i.e., ROEt).

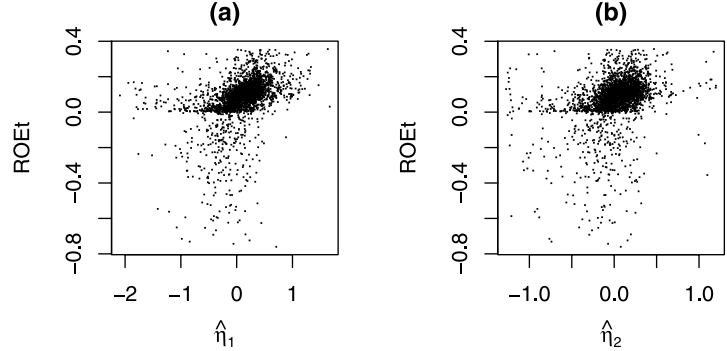

Fig. 1. *Scatter plots of ROEt versus $\hat{\eta}_j$ $(j = 1, 2)$.*

Although both CP-DR directions depict the monotonically increasing pattern when ROEt $> 0$, there is no clear pattern that can be identified for those observations with ROEt $< 0$. This is because firms operating in the Chinese stock market have extremely strong motivation to avoid reporting negative earnings. Otherwise, they might be subject to severe punishment from the China Security Regulation Commission (the government body overseeing the stock market). In addition, the firms with negative ROEt values are typically among those with relatively poor earnings capability. Thus, they are most likely to be involved in heavy earnings management. This induces a value-destroying process [Jiang and Wang (2008)] and makes the resulting earnings pattern (i.e., the regression relationship as in Figure 1) depart from the fundamental economic rules. Consequently, no clear regression pattern can be detected for those observations with ROEt $< 0$.

**8. Discussion.** In this paper, we employ the contour projected approach to establish a new theory for sufficient dimension reduction. Our approach leads to the notion of GCS, which is closely related to, but very different from, the traditional CS. To estimate the GCS, we employ three methods, CP-SIR, CP-SAVE and CP-DR, via the CP theory. Monte Carlo studies demonstrate that they are superior to SIR, SAVE and DR, respectively, especially in cases where the predictors have heavy-tailed distributions. In the development of CP theory, we mainly focus on population properties. Therefore, several important topics worth further investigation remain. The first avenue of research would be to further establish the asymptotic normality of the three contour projected estimators [Zhu and Fang (1996) and Li and Zhu (2007)] as well as to investigate the consistency of the structural dimension estimator, MERC. The second area would be to extend the contour projected approach to existing dimension reduction methods such as dimension reduction via higher-order moments, [Yin and Cook (2002, 2003, 2004)], dimension reduction in multivariate regressions, [Cook and Setodji (2003), Li et al. (2003) and Li, Wen and Zhu (2008)] and shrinkage inverse regressions [Ni, Cook and Tsai (2005)]. We believe that these efforts would enhance the usefulness of CP in dimension reduction.

## APPENDIX

### A.1. Proof of Theorem 1.

(SUFFICIENCY OF CONTOUR ASYMMETRIC). According to Proposition 6.2 of Cook (1998b), $\mathcal{C}_{y|\overrightarrow{x}}$ is unique, if it exists. Hence, we only need to establish the existence of the CCS. Furthermore, applying the results of Cook (1996) and Cook (1998b), it is equivalent to show that the intersection [denoted by $\mathcal{S}(\delta)$] of two arbitrary SCSs [denoted by $\mathcal{S}(\alpha)$ and $\mathcal{S}(\beta)$] is still a SCS.

If $\mathcal{S}(\alpha) \subset \mathcal{S}(\beta)$ or $\mathcal{S}(\alpha) \supset \mathcal{S}(\beta)$, then $\mathcal{S}(\alpha) \cap \mathcal{S}(\beta)$ is a SCS. Thus, we only need to consider $\mathcal{S}(\alpha) \not\subset \mathcal{S}(\beta)$ and $\mathcal{S}(\alpha) \not\supset \mathcal{S}(\beta)$, which indicates that $\mathcal{S}(\delta) \neq \mathcal{S}(\alpha)$ and $\mathcal{S}(\delta) \neq \mathcal{S}(\beta)$. As a result, the bases $\alpha$ and $\beta$ can be further decomposed as $\alpha = (\alpha_1, \delta)$ and $\beta = (\beta_1, \delta)$, respectively, for some $\alpha_1 \neq \varnothing$ and $\beta_1 \neq \varnothing$. Then, define $W = (W_1, W_2, W_3) = (\alpha_1^\top \overrightarrow{X}, \beta_1^\top \overrightarrow{X}, \delta^\top \overrightarrow{X}) = \eta^\top \overrightarrow{X}$, where the support of $\overrightarrow{X}$ is the unit contour $\{\overrightarrow{x} : \|\overrightarrow{x}\|^2 = 1\}$. Consequently, the support of $W$ is either a unit contour, if $\text{rank}(\eta) = p$, or the convex set $\{w : \|w\| \leq 1\}$ if $\text{rank}(\eta) < p$. This allows us to consider the two separate cases given below to prove that $\mathcal{S}(\delta)$ is a SCS.

CASE 1 $[\text{rank}(\eta) < p]$. Let $\Omega_{12|3}(w_3) = \{(w_1, w_2) : \|w_1\|^2 + \|w_2\|^2 \leq 1 - \|w_3\|^2\}$, which is the support of the conditional distribution $(W_1, W_2)|W_3 = w_3$. For any $(w_1, w_2) \in \Omega_{12|3}(w_3)$, we have $\|w_1\|^2 + 0^2 \leq \|w_1\|^2 + \|w_2\|^2 \leq 1 - \|w_3\|^2$. Hence, $(w_1, 0) \in \Omega_{12|3}(w_3)$. Because both $\mathcal{S}(\eta)$ and $\mathcal{S}(\alpha)$ are SCSs, we obtain

$$\begin{aligned}
G_y(\overrightarrow{X} = \overrightarrow{x}) &= G_y(W_1 = w_1, W_2 = w_2, W_3 = w_3) \\
&= G_y(W_1 = w_1, W_3 = w_3) \\
&= G_y(W_1 = w_1, W_2 = 0, W_3 = w_3).
\end{aligned} \quad (A.1)$$

Let $(w_1^*, w_2^*)$ be another arbitrary point in $\Omega_{12|3}$ such that $\|w_1^*\|^2 + \|w_2^*\|^2 \leq 1 - \|w_3\|^2$. This implies that $\|w_1^*\|^2 + 0^2 \leq 1 - \|w_3\|^2$, which leads to $(w_1^*, 0) \in \Omega_{12|3}(w_3)$. This, together with the fact that $\mathcal{S}(\beta)$ is a SCS, results in

$$\begin{aligned}(\text{A.2}) \quad G_y(W_1 = w_1, W_2 = 0, W_3 = w_3) &= G_y(W_2 = 0, W_3 = w_3) \\ &= G_y(W_1 = w_1^*, W_2 = 0, W_3 = w_3).\end{aligned}$$

Equations (A.1) and (A.2) yield

$$(\text{A.3}) \quad G_y(\vec{X} = \vec{x}) = G_y(W_1 = w_1^*, W_2 = 0, W_3 = w_3).$$

Because $\mathcal{S}(\alpha)$ is a SCS, we obtain

$$(\text{A.4}) \quad G_y(W_1 = w_1^*, W_2 = 0, W_3 = w_3) = G_y(W_1 = w_1^*, W_3 = w_3),$$

$$(\text{A.5}) \quad G_y(W_1 = w_1^*, W_2 = w_2^*, W_3 = w_3) = G_y(W_1 = w_1^*, W_3 = w_3).$$

By (A.1), (A.3), (A.4) and (A.5), we have

$$G_y(W_1 = w_1, W_2 = w_2, W_3 = w_3) = G_y(W_1 = w_1^*, W_2 = w_2^*, W_3 = w_3).$$

This implies that $G_y(W_1 = w_1, W_2 = w_2, W_3 = w_3)$ is a constant function of $(W_1, W_2) \in \Omega_{12|3}(w_3)$. As a result,

$$\begin{aligned}(\text{A.6}) \quad G_y(\vec{X} = \vec{x}) &= G_y(W_1 = w_1, W_2 = w_2, W_3 = w_3) \\ &= G_y(W_3 = w3) = G_y(\delta^T \vec{X} = \delta^T \vec{x}).\end{aligned}$$

Applying Lemma 1 of Zeng and Zhu (2008), (A.6) leads to the conclusion that $\mathcal{S}(\delta)$ is a SCS. Note that Cook (1996, 1998b) has used this constant function technique to prove his Lemma 2 and Proposition 6.4, respectively.

CASE 2 [rank$(\eta) = p$]. Let $\Omega_W = \{w : \|w_1\|^2 + \|w_2\|^2 + \|w_3\|^2 = 1\}$, which is a unit contour. In addition, let $\tilde{w}_j = w_j / \|w_j\|$ for $j = 1, \ldots, 3$. Because $\Omega_W$ is a unit contour, any component of $(\|W_1\|, \|W_2\|, \|W_3\|)$ is uniquely determined by the other two components. Therefore, for any $w \in \Omega_W$, we have

$$\begin{aligned}(\text{A.7}) \quad & G_y(W = w) \\ &= G_y(\tilde{W}_1 = \tilde{w}_1, \tilde{W}_2 = \tilde{w}_2, \|W_1\| = \|w_1\|, \|W_2\| = \|w_2\|, W_3 = w_3) \\ &= G_y(\tilde{W}_1 = \tilde{w}_1, \tilde{W}_2 = \tilde{w}_2, \|W_1\| = \|w_1\|, W_3 = w_3).\end{aligned}$$

Let $(w_1^*, w_2^*)$ be an another arbitrary point in $\Omega_W$ such that $\|w_1^*\|^2 + \|w_2^*\|^2 + \|w_3^*\|^2 = 1$. In addition, let $\tilde{w}_j^* = w_j^* / \|w_j^*\|$ for $j = 1, \ldots, 3$. Because $\mathcal{S}(\alpha)$ is a SCS, we are able to apply the same techniques as used in the proof of Case 1 and equation (A.7) to obtain

$$G_y(W = w) = G_y(\tilde{W}_1 = \tilde{w}_1, \tilde{W}_2 = \tilde{w}_2^*, \|W_1\| = \|w_1\|, W_3 = w_3).$$

Using the fact that $\Omega_W$ is a unit contour, the above equation can be expressed as

$$G_y(W = w)$$
$$= G_y(\tilde{W}_1 = \tilde{w}_1, \tilde{W}_2 = \tilde{w}_2^*, \|W_2\| = (1 - \|w_1\|^2 - \|w_3\|^2)^{1/2}, W_3 = w_3).$$

Moreover, $\mathcal{S}(\beta)$ is a SCS. This allows us to employ the same technique used in the proof of Case 1 to express the above equation as

$$G_y(W = w) = G_y(\tilde{W}_1 = \tilde{w}_1^*, \tilde{W}_2 = \tilde{w}_2^*,$$
$$\text{(A.8)} \qquad \|W_2\| = (1 - \|w_1\|^2 - \|w_3\|^2)^{1/2}, W_3 = w_3)$$
$$= G_y(\tilde{W}_1 = \tilde{w}_1^*, \tilde{W}_2 = \tilde{w}_2^*, \|W_1\| = \|w_1\|, W_3 = w_3).$$

Note that $(\tilde{w}_1, \tilde{w}_2)$ and $(\tilde{w}_1^*, \tilde{w}_2^*)$ are arbitrary points. Hence, applying the same argument used in the proof of Case 1, (A.7) and (A.8) lead to

$$\text{(A.9)} \qquad G_y(W = w) = G_y(\|W_1\| = \|w_1\|, W_3 = w_3).$$

Because we assume that $Y|\overrightarrow{X}$ is *contour asymmetric*, $G_y(W = w)$ must be degenerate on $\|W_1\|$. As a result, $\mathcal{S}(\delta)$ is a SCS.

(NECESSITY OF CONTOUR ASYMMETRIC). To show that *contour asymmetric* is a necessary condition for the existence of the CCS, it suffices to show that $\mathcal{C}_{y|\overrightarrow{x}}$ does not exist, as long as $Y|\overrightarrow{X}$ is *contour symmetric* on some direction. To this end, we assume that $Y|\overrightarrow{X}$ is *contour symmetric* on direction $B_1$. Then, there exists another direction $B_2$ satisfying

$$\text{(A.10)} \quad G_y(\overrightarrow{X} = \overrightarrow{x}) = G_y(\|P_{B_1}\overrightarrow{X}\| = \|P_{B_1}\overrightarrow{x}\|, P_{B_2}\overrightarrow{X} = P_{B_2}\overrightarrow{x}),$$

where $\mathcal{S}(B_1) \cap \mathcal{S}(B_2) = \varnothing$, $\mathcal{S}(B_1) \cup \mathcal{S}(B_2) \neq \mathbb{R}^p$, and $G_y(\overrightarrow{X} = \overrightarrow{x})$ is a non-degenerate function in $\|P_{B_1}\overrightarrow{x}\|$.

It is easy to show that (A.10) implies that $\mathcal{S}(\tilde{B}) = \mathcal{S}(B_1) \cup \mathcal{S}(B_2)$ is a SCS. Furthermore, because the support of $\overrightarrow{X}$ is a unit contour, the right-hand side of (A.10) can be rewritten as

$$G_y(\overrightarrow{X} = \overrightarrow{x}) = G_y(\|P_{B_{1,2}^\perp}\overrightarrow{X}\| = \|P_{B_{1,2}^\perp}\overrightarrow{x}\|, P_{B_2}\overrightarrow{X} = P_{B_2}\overrightarrow{x}),$$

where $B_{1,2}^\perp$ is a basis of the orthogonal subspace of $\mathcal{S}(B_1) \cup \mathcal{S}(B_2)$, and

$$\{\mathcal{S}(B_1) \cup \mathcal{S}(B_2)\} \cup \mathcal{S}(B_{1,2}^\perp) = \mathbb{R}^p.$$

As a result, $\mathcal{S}(\tilde{B}^*) = \mathcal{S}(B_{1,2}^\perp) \cup \mathcal{S}(B_2)$ is also a SCS of $Y|\overrightarrow{X}$. However, $\mathcal{S}(\tilde{B}^*) \cap \mathcal{S}(\tilde{B}) = \mathcal{S}(B_2)$ is not a SCS. Otherwise, by (A.10), $G_y(\overrightarrow{X} = \overrightarrow{x})$ would be degenerate on $\|P_{B_1}\overrightarrow{x}\|$, which contradicts with the definition of *contour symmetry*. This completes the proof of necessity.

**A.2. Proof of Theorem 2.** Without loss of generality, we assume that there are two different GCSs ($\mathcal{S}_1$ and $\mathcal{S}_2$) with the minimal structural dimension $d_0$ [i.e., $\dim(\mathcal{S}_1) = \dim(\mathcal{S}_2) = d_0$]. Because $\mathcal{S}_1 \neq \mathcal{S}_2$, we have $\dim(\mathcal{S}_1 \cap \mathcal{S}_2) < d_0$. Furthermore, according to the assumption that the structural dimensions of $\mathcal{S}_1$ and $\mathcal{S}_2$ are minimal, $\mathcal{S}_1 \cap \mathcal{S}_2$ cannot be a SCS. Thus, if we are able to show that $\mathcal{S}_1 \cap \mathcal{S}_2$ is also a SCS, then the first part of Theorem 2 follows. Based on the definition of $\mathcal{K}_{y|\vec{x}}$, every SCS must contain $\mathcal{K}_{y|\vec{x}}$, which implies the following inequality:

$$\dim(\mathcal{S}_1 \cup \mathcal{S}_2) \leq \dim(\mathcal{S}_1) + \dim(\mathcal{S}_2) - \dim(\mathcal{K}_{y|\vec{x}}).$$

In addition, according to the condition $d_0 < \{p + \dim(\mathcal{K}_{y|\vec{x}})\}/2$, we have

$$\dim(\mathcal{S}_1) + \dim(\mathcal{S}_2) - \dim(\mathcal{K}_{y|\vec{x}})$$
$$< 2 \times [\{p + \dim(\mathcal{K}_{y|\vec{x}})\}/2] - \dim(\mathcal{K}_{y|\vec{x}}) = p.$$

The above two equations together imply that $\dim(\mathcal{S}_1 \cup \mathcal{S}_2) < p$. Moreover, applying the same techniques used in the proof of Case 1 in Theorem 1, we obtain the result that $\mathcal{S}_1 \cap \mathcal{S}_2$ is a SCS. The proof of the first part is complete.

We next show the second part of Theorem 2. Let $\mathcal{G}_{y|\vec{x}}$ be an arbitrary GCS. Because the CCS $\mathcal{C}_{y|\vec{x}}$ exists, we have $\mathcal{G}_{y|\vec{x}} \supset \mathcal{C}_{y|\vec{x}}$. On the other hand, the GCS is the SCS with the minimal structural dimension and the CCS is also a SCS. Accordingly, $\dim(\mathcal{G}_{y|\vec{x}}) \leq \dim(\mathcal{C}_{y|\vec{x}})$, which leads to $\mathcal{G}_{y|\vec{x}} = \mathcal{C}_{y|\vec{x}}$. This, together with the uniqueness of $\mathcal{C}_{y|\vec{x}}$, implies $\mathcal{G}_{y|\vec{x}}$ is unique. The proof of the second part is complete.

**A.3. Proof of Theorem 3.**

STATEMENT (1). Let $\mathcal{B} = \{B : Y \perp\!\!\!\perp X | B^\top X\}$ and $\mathcal{B}^* = \{B : Y \perp\!\!\!\perp \vec{X} | B^\top \vec{X}\}$. Applying the Lemma 2 of Wang, Ni and Tsai (2008), we have $\mathcal{B} \subset \mathcal{B}^*$. Accordingly,

$$\mathcal{K}_{y|\vec{x}} = \bigcap_{B \in \mathcal{B}^*} \mathcal{S}(B) = \left\{ \bigcap_{B \in \mathcal{B}} \mathcal{S}(B) \right\} \cap \left\{ \bigcap_{B \in \mathcal{B}^* \setminus \mathcal{B}} \mathcal{S}(B) \right\}$$
$$= \mathcal{S}_{y|x} \cap \left\{ \bigcap_{B \in \mathcal{B}^* \setminus \mathcal{B}} \mathcal{S}(B) \right\} \subset \mathcal{S}_{y|x}.$$

This completes the proof.

STATEMENT (2). By Lemma 2 of Wang, Ni and Tsai (2008), we know that $\mathcal{S}_{y|x}$ is also a SCS. Then, by the definition of the GCS, we have that the structural dimension of $\mathcal{G}_{y|\vec{x}}$ cannot be larger than that of $\mathcal{S}_{y|x}$. The proof is complete.

**A.4. Proof of Theorem 4.** To prove the theorem, we consider two different cases, namely $Y|\overrightarrow{X}$ is *contour asymmetric* and $\dim(\mathcal{S}_{y|x}) < \{p + \dim(\mathcal{K}_{y|\overrightarrow{x}})\}/2$.

CASE 1 ($Y|\overrightarrow{X}$ is contour asymmetric). Applying Theorem 1 and the second part of Theorem 2, we have $\mathcal{G}_{y|\overrightarrow{x}} = \mathcal{C}_{y|\overrightarrow{x}} = \mathcal{K}_{y|\overrightarrow{x}}$. To obtain $\mathcal{G}_{y|\overrightarrow{x}} = \mathcal{S}_{y|x}$, it suffices to show that $\mathcal{K}_{y|\overrightarrow{x}} = \mathcal{S}_{y|x}$. If $\mathcal{K}_{y|\overrightarrow{x}} = \mathcal{C}_{y|\overrightarrow{x}} = \varnothing$, then $Y \perp\!\!\!\perp \overrightarrow{X}$. This implies that $Y = g(R, \varepsilon)$, where $g(\cdot)$ is an unknown function and $\varepsilon$ is independent of $\overrightarrow{X}$ and $R$. Accordingly, $\mathcal{S}_{y|x} = \varnothing$ or $\mathcal{S}_{y|x} = \mathbb{R}^p$. However, the latter situation contradicts the condition that $Y|X$ is *dimension reducible*. Thus, we only consider the situation where $\dim(\mathcal{K}_{y|\overrightarrow{x}}) > 0$.

Let $B$ be a basis of $\mathcal{K}_{y|\overrightarrow{x}} = \mathcal{C}_{y|\overrightarrow{x}}$. Based on Theorem 3(1), it suffices to show that $\mathcal{S}(B)$ is a SDR subspace of $Y|X$. From $Y \perp\!\!\!\perp \overrightarrow{X}|B^\top \overrightarrow{X}$, we have $Y \perp\!\!\!\perp (\overrightarrow{X}, R)|(P_B \overrightarrow{X}, R)$, which is equivalent to

$$G_y(\overrightarrow{X} = \overrightarrow{x}, R = r) = G_y(P_B \overrightarrow{X} = P_B \overrightarrow{x}, R = r).$$

Note that $(P_B \overrightarrow{X}, R)$ is a one-to-one mapping of $(P_B X, \|P_{B^\perp} X\|)$, where $B^\perp$ is a basis satisfying $\mathcal{S}(B) \cup \mathcal{S}(B^\perp) = \mathbb{R}^p$ and $\mathcal{S}(B) \cap \mathcal{S}(B^\perp) = \varnothing$. In addition, $(\overrightarrow{X}, R)$ a one-to-one mapping of $X$. As a result, we obtain

$$G_y(X = x) = G_y(\overrightarrow{X} = \overrightarrow{x}, R = r)$$
$$= G_y(P_B \overrightarrow{X} = P_B \overrightarrow{x}, R = r)$$
$$= G_y(P_B X = P_B x, \|P_{B^\perp} X\| = \|P_{B^\perp} x\|),$$

where $x = r\overrightarrow{x}$. Then, Theorem 4 follows if we are able to show that $G_y(X = x)$ is degenerate on $\|P_{B^\perp} x\|$.

We assume that $Y|X$ is *dimension reducible*. Hence, if $G_y(X = x)$ is not degenerate on $\|P_{B^\perp} x\|$, then $G_y(X = x)$ must be degenerate on $\mathcal{S}(B_1)$, which is a subspace of $\mathcal{S}(B)$ with $\dim\{\mathcal{S}(B_1)\} > 0$. Accordingly,

$$G_y(X = x) = G_y(P_{B_2} X = P_{B_2} x, \|P_{B^\perp} X\| = \|P_{B^\perp} x\|),$$

where $B_2$ satisfies $\mathcal{S}(B_1) \cup \mathcal{S}(B_2) = \mathcal{S}(B)$ and $\mathcal{S}(B_1) \cap \mathcal{S}(B_2) = \varnothing$. This implies that $\mathcal{S}(B^*) = \{\mathcal{S}(B_2) \cup \mathcal{S}(B^\perp)\}$ is a SDR subspace of $Y|X$. Then, from Lemma 2 of Wang, Ni and Tsai (2008), $\mathcal{S}(B^*)$ is also a SCS of $Y|\overrightarrow{X}$. Because $\mathcal{S}(B_1) \subset \mathcal{S}(B)$ and $\dim\{\mathcal{S}(B_1)\} > 0$, we have $\mathcal{K}_{y|\overrightarrow{x}} = \mathcal{S}(B) \not\subset \mathcal{S}(B^*)$, which contradicts the definition of $\mathcal{K}_{y|\overrightarrow{x}}$. Consequently, $G_y(X = x)$ is degenerate on $\|P_{B^\perp} x\|$ and $\mathcal{K}_{y|\overrightarrow{x}} = \mathcal{S}_{y|x}$.

CASE 2 [$\dim(\mathcal{S}_{y|x}) < \{p + \dim(\mathcal{K}_{y|\overrightarrow{x}})\}/2$]. Similar to Case 1, we only consider a situation where $\dim(\mathcal{G}_{y|\overrightarrow{x}}) > 0$. By Lemma 2 of Wang, Ni and Tsai (2008), we know that $\mathcal{S}_{y|x}$ is also a SCS. Then, according to the definition of GCS, we have $\dim(\mathcal{G}_{y|\overrightarrow{x}}) \leq \dim(\mathcal{S}_{y|x})$. As a result, we obtain $\mathcal{G}_{y|\overrightarrow{x}} \subset$ $\mathcal{S}_{y|x}$. Otherwise, applying the condition that $\dim(\mathcal{G}_{y|\vec{x}}) \leq \dim(\mathcal{S}_{y|x}) < \{p + \dim(\mathcal{K}_{y|\vec{x}})\}/2$, together with the same technique used in the proof of Theorem 2, we can show that $\{\mathcal{G}_{y|\vec{x}} \cap \mathcal{S}_{y|x}\} \neq \mathcal{G}_{y|\vec{x}}$, but is a SCS. This indicates that the structural dimension of $\{\mathcal{G}_{y|\vec{x}} \cap \mathcal{S}_{y|x}\}$ is smaller than that of $\mathcal{G}_{y|\vec{x}}$, which contradicts the definition of GCS. Thus, we must have $\mathcal{G}_{y|\vec{x}} \subset \mathcal{S}_{y|x}$.

We slightly abuse notation by letting $B$ be a basis of $\mathcal{G}_{y|\vec{x}}$, which is different from Case 1. Following the same argument used in Case 1, we have $G_y(X = x) = G_y(P_B X = P_B x, \|P_{B^\perp} X\| = \|P_{B^\perp} x\|)$. Recall that $Y|X$ is *dimension reducible*. Hence, if $G_y(X = x)$ is not degenerate on $\|P_{B^\perp} x\|$, then there must exist two bases $B_1$ and $B_2$ such that $G_y(X = x)$ is degenerate on $P_{B_2} x$, $\text{rank}(B_2) > 0$, $\mathcal{S}(B_1) \cup \mathcal{S}(B_2) = \mathcal{S}(B)$, and $\mathcal{S}(B_1) \cap \mathcal{S}(B_2) = \varnothing$. As a result, $G_y(X = x) = G_y(P_{B_1} X = P_{B_1} x, \|P_{B^\perp} X\| = \|P_{B^\perp} x\|)$, which implies that $\mathcal{S}^* = \mathcal{S}(B_1) \cup \mathcal{S}(B^\perp)$ is a SDR subspace. Moreover, $\mathcal{S}(B_2) \not\subset \mathcal{S}^*$, and hence $\mathcal{S}(B) = \mathcal{G}_{y|\vec{x}} \not\subset \mathcal{S}^*$. However, according to the definition of CS, we must have $\mathcal{S}_{y|x} \subset \mathcal{S}^*$, and hence $\mathcal{G}_{y|x} \subset \mathcal{S}^*$ by combining the fact $\mathcal{G}_{y|\vec{x}} \subset \mathcal{S}_{y|x}$. Therefore, we will get a contradiction if $G_y(X = x)$ is degenerate on $P_{B_2} x$. Consequently, $G_y(X = x)$ is degenerate on $\|P_{B^\perp} x\|$ which implies $G_y(X = x) = G_y(P_B X = P_B x)$, which implies that $\mathcal{G}_{y|\vec{x}}$ is also a SDR subspace. Combing the previous result $\mathcal{G}_{y|\vec{x}} \subset \mathcal{S}_{y|x}$, we have $\mathcal{G}_{y|\vec{x}} = \mathcal{S}_{y|x}$. The results of Cases 1 and 2 complete the proof.

**A.5. Proof of Lemma 1.** Let $\mathcal{S}(B)$ be an arbitrary SCS. Applying Lemma 1 of Wang, Ni and Tsai (2008), we have $E(\vec{X}|B^\top \vec{X}) = P_B \vec{X}$, which shows that the *linearity condition* [Li (1991)] holds on the contour projected predictor $\vec{X}$. As a result,

$$E(\vec{X}|Y) = E[E(\vec{X}|B^\top \vec{X}, Y)|Y] = E[E(\vec{X}|B^\top \vec{X})|Y] = E(P_B \vec{X}|Y) \in \mathcal{S}(B).$$

Because $\mathcal{S}(B)$ is an arbitrary SCS, we immediately have $\mathcal{S}\{E(\vec{X}|Y)\} \subset \mathcal{K}_{y|\vec{x}}$. This completes the proof.

**A.6. Proof of Lemma 2.**

STATEMENT (1). For any $\mathcal{S}(B) \in \mathcal{I}$, if $Y|\vec{X}$ is *contour asymmetric*, then the desired result follows because $\mathcal{S}(B_A) = \mathcal{S}(B)$. In contrast, if $Y|\vec{X}$ is *contour symmetric*, then $\mathcal{S}(B_S) \cap \mathcal{K}_{y|\vec{x}} = \varnothing$ (see the proof of *necessity of contour asymmetric* in Appendix A.1). This implies that $\mathcal{K}_{y|\vec{x}} = \bigcap_{\mathcal{S}(B) \in \mathcal{I}} \mathcal{S}(B_A)$.

STATEMENT (2). For any $\mathcal{S}(B) \in \mathcal{I}$, if $Y|\vec{X}$ is *contour symmetric* on some direction $B_S$, then $\mathcal{K}_{y|\vec{x}} \subset \mathcal{S}(B_A)$ (see the proof of *necessity of contour asymmetric* in Appendix A.1). It implies that $\mathcal{K}_{y|\vec{x}} \neq \mathcal{G}_{y|\vec{x}}$. Otherwise, $\mathcal{K}_{y|\vec{x}}$ should be a SCS. Then, by previous statement, we know that $\mathcal{S}(B_A) \supset$

$\mathcal{K}_{y|\overrightarrow{x}}$. Hence $\mathcal{S}(B_A)$ is also a SCS. Consequently, $G_y(\overrightarrow{X} = \overrightarrow{x})$ degenerates on $P_{B_S}\overrightarrow{x}$, which contradicts the assumption that $Y|\overrightarrow{X}$ is *contour symmetric* on $B_S$.

STATEMENT (3). If $Y|\overrightarrow{X}$ is *contour asymmetric*, then Theorem 1 and the second part of Theorem 2 lead to $\mathcal{K}_{y|\overrightarrow{x}} = \mathcal{C}_{y|\overrightarrow{x}} = \mathcal{G}_{y|\overrightarrow{x}}$. The proof is complete.

**A.7. Proof of Theorem 5.** Let $M_{\text{SIR}} = \text{cov}\{(E(\overrightarrow{X}|Y)\}$, which is the kernel matrix of CP-SIR. By Lemmas 1 and 2, we have $\mathcal{S}(M_{\text{SIR}}) \subset \mathcal{K}_{y|\overrightarrow{x}} = \mathcal{C}_{y|\overrightarrow{x}} = \mathcal{G}_{y|\overrightarrow{x}}$. Next, applying the technique from the proof of Theorem 3 in Li and Wang (2007), Theorem 5 follows if we are able to show that $v^\top M_{\text{SIR}} v > 0$ for all $v \in \mathcal{G}_{y|\overrightarrow{x}}$ with $\|v\| = 1$. Because $E(\overrightarrow{X}) = 0$, we obtain

$$v^\top M_{\text{SIR}} v = v^\top E[E(\overrightarrow{X}|Y)E(\overrightarrow{X}^\top|Y)]v = E\{[E(v^\top \overrightarrow{X}|Y)]^2\}.$$

By Assumption 1, $E(v^\top \overrightarrow{X}|Y)$ is a nondegenerate function in $Y$. In conjunction with Jensen's inequality, this leads to

$$E\{[E(v^\top \overrightarrow{X}|Y)]^2\} > \{E[E(v^\top \overrightarrow{X}|Y)]\}^2 = \{v^\top E(\overrightarrow{X})\}^2 = 0.$$

Accordingly, $v^\top M_{\text{SIR}} v > 0$ for all $v \in \mathcal{G}_{y|\overrightarrow{x}}$. This completes the proof.

**A.8. Proof of Lemma 3.** Using the fact that $E(\overrightarrow{X}|B_0\overrightarrow{X}) = P_{B_0}\overrightarrow{X}$, and given that $\mathcal{S}(B_0)$ is the GCS, we have

$$\text{cov}(\overrightarrow{X}|Y) = E[\text{cov}(\overrightarrow{X}|B_0^\top \overrightarrow{X})|Y] + P_{B_0}\text{cov}(\overrightarrow{X}|Y)P_{B_0}.$$

Let $B_0^\perp$ denote a orthonormal basis of the subspace that is the orthogonal complement of $\mathcal{S}(B_0)$. Furthermore, let $C = (B_0, B_0^\perp)$ so that $CC^\top = I_p$. Moreover, define

$$W = C^\top \overrightarrow{X} = \begin{pmatrix} B_0^\top \overrightarrow{X} \\ B_0^{\perp\top} \overrightarrow{X} \end{pmatrix} = \begin{pmatrix} W_0 \\ W_0^\perp \end{pmatrix}.$$

As a result,

$$\text{cov}(\overrightarrow{X}|B_0^\top \overrightarrow{X}) = \text{cov}(CC^\top \overrightarrow{X}|B_0^\top \overrightarrow{X}) = C \times \text{cov}(W|W_0) \times C^\top$$
$$\text{(A.11)} \qquad = C \times \begin{bmatrix} 0 & 0 \\ 0 & \text{cov}(W_0^\perp|W_0) \end{bmatrix} \times C^\top.$$

By noting that $\|W\|^2 = \|\overrightarrow{X}\|^2 = 1$, $\|W_0^\perp\|^2 = 1 - \|W_0\|^2 = 1 - \|B_0^\top \overrightarrow{X}\|^2$, and $\text{cov}(W) = \text{cov}(\overrightarrow{X}) = p^{-1}I_p$ [see Wang, Ni and Tsai (2008)], we obtain

$$\text{(A.12)} \qquad \text{cov}(W|W_0) = \frac{1 - \|W_0\|^2}{p - d_0} I_{p - d_0} = \frac{1 - \|B_0^\top \overrightarrow{X}\|^2}{p - d_0} I_{p - d_0}.$$

Applying (A.11) and (A.12), we have

$$\operatorname{cov}(\overrightarrow{X}|Y) = C \times \begin{bmatrix} 0 & 0 \\ 0 & \dfrac{1-E(\|B_0^\top \overrightarrow{X}\|^2|Y)}{p-d_0} I_{p-d_0} \end{bmatrix} \times C^\top + P_{B_0} \operatorname{cov}(\overrightarrow{X}|Y) P_{B_0}$$

$$= \frac{1-\lambda(Y)}{p-d_0} B_0^\perp (B_0^\perp)^\top + P_{B_0} \operatorname{cov}(\overrightarrow{X}|Y) P_{B_0}$$

$$= \tau(Y) Q_{B_0} + P_{B_0} \operatorname{cov}(\overrightarrow{X}|Y) P_{B_0}.$$

This implies that $\{\tau(Y)I_p - \operatorname{cov}(\overrightarrow{X}|Y)\} = P_{B_0}\{\tau(Y)I_p - \operatorname{cov}(\overrightarrow{X}|Y)\}P_{B_0}$. Subsequently, together with Lemma 1 and Ye and Weiss (2003), Lemma 3, this yields

$$\mathcal{S}\{\tau(Y)I_p - E(\overrightarrow{X}\overrightarrow{X}^\top|Y)\}$$
$$= \mathcal{S}\{\tau(Y)I_p - \operatorname{cov}(\overrightarrow{X}|Y) + E(\overrightarrow{X}|Y)E(\overrightarrow{X}^\top|Y)\} \subset \mathcal{G}_{y|\vec{x}}.$$

This completes the proof.

**A.9. Proof of Theorem 6.** Let $H = \tau(Y)I_p - E(\overrightarrow{X}\overrightarrow{X}^\top|Y)$, then the kernel matrix of CP-SAVE is $M_{\text{SAVE}} = E(H^2)$. Because $E[\tau(Y)] = 1/p$, $E(\overrightarrow{X}) = 0$, and $\operatorname{cov}(\overrightarrow{X}) = I_p/p$, we have $E(H) = 0$. Applying Lemma 3 and Ye and Weiss (2003), Lemma 3, we obtain $\mathcal{S}(M_{\text{SAVE}}) \subset \mathcal{G}_{y|\vec{x}}$. Next, if we are able to show that $v^\top M_{\text{SAVE}} v > 0$ for all $v \in \mathcal{G}_{y|\vec{x}}$ with $\|v\| = 1$, then Theorem 6 follows (the same technique has been used in the proof of Theorem 5). It is easy to see that $v^\top M_{\text{SAVE}} v = v^\top E[H(I_p - vv^\top)H]v + E[(v^\top Hv)^2]$. Because $I_p - vv^\top$ in a nonnegative definite matrix, the first term in the right-hand side of the above equation is nonnegative. By Jensen's inequality and $E(H) = 0$,

(A.13) $$E[(v^\top Hv)^2] > [E(v^\top Hv)]^2 = 0.$$

The strict inequality in (A.13) holds because $v^\top Hv$ is nondegenerate, which is shown below. After algebraic simplification, we have

$$v^\top Hv = v^\top[\tau(Y)I_p - E(\overrightarrow{X}\overrightarrow{X}^\top|Y)]v$$
$$= \tau(Y) - E[(v^\top \overrightarrow{X})^2|Y] = \tau(Y) - \phi(Y, v).$$

Note that $\tau(Y) = \{1-\lambda(Y)\}/(p-d_0)$ and

$$\lambda(Y) = E(\overrightarrow{X}^\top B_0 B_0^\top \overrightarrow{X}|Y) = \sum_{i=1}^{d_0} E(\beta_{0i}^\top \overrightarrow{X}\overrightarrow{X}^\top \beta_{0i}|Y) = \sum_{i=1}^{d_0} \phi(Y, \beta_{0i}).$$

Without loss of generality, we assume that $v = \beta_{01}$. Otherwise, we can construct a basis $B_0$ whose first column is $v$. As a result,

$$v^\top Hv = \frac{1 - \sum_{i=1}^{d_0} \phi(Y, \beta_{0i})}{p - d_0} - \phi(Y, \beta_{01})$$

$$= -(p-d_0)^{-1}\left[(p-d_0+1)\phi(Y,\beta_{01}) + \sum_{i=2}^{d_0}\phi(Y,\beta_{0i}) - 1\right],$$

which is a linear combination of $\phi(Y,\beta_{0i})$, $i=1,\ldots,d_0$. According to Assumption 3, $v^\top H v$ is nondegenerate. This completes the proof.

**A.10. Proof of Theorem 7.** Because $(Y^*,\overrightarrow{X}^*)$ is an independent copy of $(Y,\overrightarrow{X})$, we have

$$\begin{aligned}
A(Y,Y^*) &= E[\overrightarrow{X}\overrightarrow{X}^\top - \overrightarrow{X}(\overrightarrow{X}^*)^\top - \overrightarrow{X}^*\overrightarrow{X}^\top + \overrightarrow{X}^*(\overrightarrow{X}^*)^\top | Y,Y^*] \\
&= E[\overrightarrow{X}\overrightarrow{X}^\top | Y] - E[\overrightarrow{X}|Y]E[(\overrightarrow{X}^*)^\top | Y^*] \\
&\quad - E[\overrightarrow{X}^*|Y^*]E[\overrightarrow{X}^\top | Y] + E[\overrightarrow{X}^*(\overrightarrow{X}^*)^\top | Y^*].
\end{aligned} \tag{A.14}$$

Using (A.14), we further obtain

$$\begin{aligned}
&[\tau(Y) + \tau(Y^*)]I_p - A(Y,Y^*) \\
&= \{\tau(Y)I_p - E(\overrightarrow{X}\overrightarrow{X}^\top | Y)\} + \{\tau(Y^*)I_p - E[\overrightarrow{X}^*(\overrightarrow{X}^*)^\top | Y^*]\} \\
&\quad + \{E[\overrightarrow{X}|Y]E[(\overrightarrow{X}^*)^\top | Y^*] + E[\overrightarrow{X}^*|Y^*]E[\overrightarrow{X}^\top | Y]\}.
\end{aligned} \tag{A.15}$$

Equation (A.15), in conjunction with Lemmas 1 and 3, yields the desired result.

**A.11. Proof of Theorem 8.** Let $D = [\tau(Y) + \tau(Y^*)]I_p - A(Y,Y^*)$ so that $M_{\text{DR}} = E(D^2)$. Note that $E[\tau(Y)] = 1/p$, $E(\overrightarrow{X}) = 0$, $\text{cov}(\overrightarrow{X}) = I_p/p$, and $(\overrightarrow{X}^*,Y^*)$ is an independent copy of $(\overrightarrow{X},Y)$. Then applying the results of (A.14) and (A.15), we are able to show that $E[A(Y,Y^*)] = 2I_p/p$ and $E(D) = 0$. Furthermore, employing Theorem 7 together with Ye and Weiss (2003), Lemma 3, we have $\mathcal{S}(M_{\text{DR}}) \subset \mathcal{G}_{y|\overrightarrow{x}}$. Accordingly, the theorem follows if we can show that $v^\top M_{\text{DR}} v > 0$ for all $v \in \mathcal{G}_{y|\overrightarrow{x}}$ with $\|v\| = 1$. It is easy to see that

$$v^\top M_{\text{DR}} v = v^\top E[D(I_p - vv^\top)D]v + E[(v^\top Dv)^2].$$

Because $I_p - vv^\top$ is nonnegative definite matrix, the first term in the right-hand side of the above equation is nonnegative. By Jensen's inequality and $E(D) = 0$, we have

$$E[(v^\top Dv)^2] = \text{var}(v^\top Dv) > [E(v^\top Dv)]^2 = 0. \tag{A.16}$$

The strict inequality in (A.16) holds because $v^\top Dv$ is nondegenerate, which is shown next. Let $\varphi(Y,v) = E(v^\top \overrightarrow{X} | Y)$. Then, (A.15) implies that

$$v^\top Dv = \tau(Y) + \tau(Y^*) - \phi(Y,v) - \phi(Y^*,v) + 2\varphi(Y,v)\varphi(Y^*,v).$$

As noted in the proof of Theorem 8, we have $\tau(Y) = \{1 - \lambda(Y)\}/(p - d_0)$ and $\lambda(Y) = E(\overrightarrow{X}^\top B_0 B_0^\top \overrightarrow{X} | Y) = \sum_{i=1}^{d_0} E(\beta_{0i}^\top \overrightarrow{X} \overrightarrow{X}^\top \beta_{0i} | Y) = \sum_{i=1}^{d_0} \phi(Y, \beta_{0i})$. Without loss of generality, we also assume that $v = \beta_{01}$. Accordingly,

$$v^\top Dv = \frac{1 - \sum_{i=1}^{d_0} \phi(Y, \beta_{0i})}{p - d_0} + \frac{1 - \sum_{i=1}^{d_0} \phi(Y^*, \beta_{0i})}{p - d_0}$$
$$- \phi(Y, \beta_{01}) - \phi(Y^*, \beta_{01}) + 2\varphi(Y, \beta_{01})\varphi(Y^*, \beta_{01})$$

$$= -(p - d_0)^{-1}\left\{\left[(p - d_0 + 1)[\phi(Y, \beta_{01}) - p^{-1}]\right.\right.$$
$$\left.+ \sum_{i=2}^{d_0}[\phi(Y, \beta_{0i}) - p^{-1}]\right]$$
$$+ \left[(p - d_0 + 1)[\phi(Y^*, \beta_{01}) - p^{-1}]\right.$$
$$\left.\left.+ \sum_{i=2}^{d_0}[\phi(Y^*, \beta_{0i}) - p^{-1}]\right]\right\}$$
$$+ 2\varphi(Y, \beta_{01})\varphi(Y^*, \beta_{01})$$
$$= -(p - d_0)^{-1}\{H(Y) + H(Y^*)\} + 2\varphi(Y, \beta_{01})\varphi(Y^*, \beta_{01}),$$

where $H(Y) = (p - d_0 + 1)[\phi(Y, \beta_{01}) - p^{-1}] + \sum_{i=2}^{d_0}[\phi(Y, \beta_{0i}) - p^{-1}]$. It is easy to see that $E\{\phi(Y, \beta_{0i})\} = 1/p$, $E\{\varphi(Y, \beta_i)\} = 0$, and $E\{H(Y)\} = 0$. In addition,

$$\text{cov}[H(Y), \varphi(Y, v)\varphi(Y^*, v)] = E[H(Y)\varphi(Y, v)\varphi(Y^*, v)]$$
$$= E[H(Y)\varphi(Y, v)]E[\varphi(Y^*, v)] = 0.$$

Hence, we have

$$\text{var}(v^\top Dv) = \frac{\text{var}[H(Y)] + \text{var}[H(Y^*)]}{(p - d_0)^2}$$
$$(A.17) \qquad + 4\text{var}[\varphi(Y, \beta_{01})]\text{var}[\varphi(Y^*, \beta_{01})]$$
$$= \frac{2\text{var}[H(Y)]}{(p - d_0)^2} + 4\{\text{var}[\varphi(Y, \beta_{01})]\}^2.$$

If Assumption 1 holds and $Y|\overrightarrow{X}$ is *contour asymmetric*, we then have $\varphi(Y, \beta_{01})$ nondegenerate. Thus, the second term in (A.17) is strictly positive. Otherwise, if Assumption 2 holds and $Y|\overrightarrow{X}$ is *dimension reducible*, we then have $H(Y)$ nondegenerate. Therefore, the first term in (A.17) is strictly positive. As a result, $\text{var}(v^\top Dv) > 0$, or equivalently, $v^\top Dv$ is nondegenerate. This completes the proof.

**A.12. Proof of Theorem 9.** Because $(Y^*, \vec{X}^*)$ is an independent copy of $(Y, \vec{X})$ and $E[\tau(Y)] = 1/p$, we have

$$
\begin{aligned}
M_{\mathrm{DR}} &= E\{[[\tau(Y) + \tau(Y^*)]I_p - A(Y, Y^*)]^2\} \\
&= \{2E\{[\tau(Y)]^2\} + 2p^{-2}\}I_p \\
&\quad - 4E\{\tau(Y)A(Y, Y^*)\} + E\{[A(Y, Y^*)]^2\}.
\end{aligned}
\tag{A.18}
$$

Using the result that $E(\vec{X}) = 0$ and $\mathrm{cov}(\vec{X}) = I_p/p$, we can simplify $E\{\tau(Y) \times A(Y, Y^*)\}$ in the above equation as follows:

$$
\begin{aligned}
&E[\tau(Y)A(Y, Y^*)] \\
&= E[\tau(Y)E(\vec{X}\vec{X}^\top|Y)] + E[\tau(Y)]E[\vec{X}^*(\vec{X}^*)^\top] \\
&\quad - E[\tau(Y)E(\vec{X}|Y)]E[(\vec{X}^*)^\top] - E(\vec{X}^*)E[\tau(Y)E(\vec{X}^\top|Y)] \\
&= E[\tau(Y)E(\vec{X}\vec{X}^\top|Y)] + p^{-2}I_p.
\end{aligned}
\tag{A.19}
$$

Furthermore, the third term of (A.18) can be reduced to

$$
\begin{aligned}
&E\{[A(Y, Y^*)]^2\} \\
&= 2\{E[[E(\vec{X}\vec{X}^\top|Y)]^2] + [E[E(\vec{X}|Y)E(\vec{X}^\top|Y)]]^2 \\
&\quad + E[E(\vec{X}^\top|Y)E(\vec{X}|Y)]E[E(\vec{X}|Y)E(\vec{X}^\top|Y)] + p^{-2}I_p\}.
\end{aligned}
\tag{A.20}
$$

Substituting (A.20) and (A.19) into (A.18), we obtain

$$
\begin{aligned}
M_{\mathrm{DR}} &= 2\{E[\tau^2(Y)]I_p + E[E^2(\vec{X}\vec{X}^\top|Y)] + E^2[E(\vec{X}|Y)E(\vec{X}^\top|Y)] \\
&\quad + E[E(\vec{X}^\top|Y)E(\vec{X}|Y)]E[E(\vec{X}|Y)E(\vec{X}^\top|Y)] \\
&\quad - 2E[\tau(Y)E(\vec{X}\vec{X}^\top|Y)]\}.
\end{aligned}
$$

This completes the proof.

### A.13. Proof of Theorem 10.

STATEMENT (1). According to Assumption 3, we are able to find an $a > 0$, such that $h(t)t^{-\alpha} < 2C_\alpha$ for $0 \leq t < a$. Thus, we have

$$
\begin{aligned}
E\|X\|^{-4} &= \int_0^\infty h(t)t^{-2}\,dt \\
&= \int_a^\infty h(t)t^{-2}\,dt + \int_0^a h(t)t^{-\alpha} \times t^{-(2-\alpha)}\,dt \\
&\leq \int_a^\infty h(t)a^{-2}\,dt + \int_0^a 2C_\alpha t^{-(2-\alpha)}\,dt \\
&\leq a^{-2}\int_0^\infty h(t)\,dt + 2C_\alpha \int_0^a t^{-(2-\alpha)}\,dt.
\end{aligned}
\tag{A.21}
$$

Because $h(t)$ is a probability density function, the first term on the right-hand side of (A.21) is finite. By Assumption 3, we have $\alpha > 1$. Thus, the second term on the right-hand side of (A.21) is also finite. Consequently, we have $E\|X\|^{-4} < \infty$. This completes the proof.

STATEMENT (2). By the definition of $\|\cdot\|_\Sigma$-norm, one can easily verify that

$$\max_{1 \leq i \leq n} \left| \frac{\|x_i - \tilde{\mu}\|_{\tilde{\Sigma}}^2}{\|x_i\|_{\tilde{\Sigma}}^2} - 1 \right|$$

$$= \max_{1 \leq i \leq n} \left| \frac{\|\tilde{\mu}\|_{\tilde{\Sigma}}^2 - 2x_i^\top \tilde{\Sigma}^{-1} \tilde{\mu}}{\|x_i\|_{\tilde{\Sigma}}^2} \right|$$

$$\leq \max_{1 \leq i \leq n} \frac{\|\tilde{\mu}\|_{\tilde{\Sigma}}^2 + 2|x_i^\top \tilde{\Sigma}^{-1} \tilde{\mu}|}{\|x_i\|_{\tilde{\Sigma}}^2} \leq \max_{1 \leq i \leq n} \frac{\|\tilde{\mu}\|_{\tilde{\Sigma}}^2 + 2\|x_i\|_{\tilde{\Sigma}} \times \|\tilde{\mu}\|_{\tilde{\Sigma}}}{\|x_i\|_{\tilde{\Sigma}}^2}$$

$$\leq \|\tilde{\Sigma}^{-1}\| \times \|\tilde{\Sigma}\| \left\{ \max_{1 \leq i \leq n} \frac{\|\tilde{\mu}\|^2}{\|x_i\|^2} \right\} + 2(\|\tilde{\Sigma}^{-1}\| \times \|\tilde{\Sigma}\|)^{1/2} \left\{ \max_{1 \leq i \leq n} \frac{\|\tilde{\mu}\|}{\|x_i\|} \right\}$$

$$= \|\tilde{\Sigma}^{-1}\| \times \|\tilde{\Sigma}\| \left\{ \frac{\|\tilde{\mu}\|^2}{\min_i \|x_i\|^2} \right\} + 2(\|\tilde{\Sigma}^{-1}\| \times \|\tilde{\Sigma}\|)^{1/2} \left\{ \frac{\|\tilde{\mu}\|}{\min_i \|x_i\|} \right\}.$$

Therefore, Statement (2) follows if we are able to show that $\|\tilde{\mu}\|^2 / \{\min_i \|x_i\|^2\} \to_p 0$. To this end, we compute the following quantity:

$$P\left(\min_{1 \leq i \leq n} \|x_i\|^2 > cn^{-1/(\alpha+1)}\right) = P^n(\|x_i\|^2 > cn^{-1/(\alpha+1)})$$

$$= \{1 - P(\|x_i\|^2 < cn^{-1/(\alpha+1)})\}^n$$

(A.22)

$$= \left(1 - \int_0^{cn^{-1/(\alpha+1)}} h(t)\,dt\right)^n$$

$$= \left(1 - \int_0^{cn^{-1/(\alpha+1)}} t^{-\alpha} h(t) \times t^\alpha \,dt\right)^n,$$

where $c > 0$ is an arbitrary constant. By Assumption 3, we have $t^{-\alpha} h(t) > C_\alpha/2$ when $n$ is large enough. As a result, the right-hand side of (A.22) is bounded by

(A.23)
$$\left(1 - \frac{C_\alpha}{2} \int_0^{cn^{-1/(\alpha+1)}} t^\alpha \,dt\right)^n = \left(1 - \frac{C_\alpha c^{(\alpha+1)}}{2(\alpha+1)} n^{-1}\right)^n$$

$$\to \exp\left\{-\frac{C_\alpha c^{(\alpha+1)}}{2(\alpha+1)}\right\}.$$

Moreover, for any arbitrarily small $\epsilon > 0$, we can find a sufficiently large $c$ such that the right-hand side of (A.23) is smaller than $\epsilon/2$. Hence, we have that

$$\limsup_{n \to \infty} P\Big(\min_{1 \leq i \leq n} \|x_i\|^2 > cn^{-1/(\alpha+1)}\Big) < \epsilon,$$

which implies $\min_{1 \leq i \leq n} \|x_i\|^2 = O_p(n^{-1/(\alpha+1)})$. This together with the fact $\|\tilde{\mu}\| = O_p(n^{-1/2})$ leads to

$$\frac{\|\tilde{\mu}\|^2}{\min_{1 \leq i \leq n} \|x_i\|^2} = O_p(n^{-1} \times n^{1/(\alpha+1)}) = O_p(n^{-\alpha/(\alpha+1)}) = o_p(1).$$

The proof is complete.

**A.14. Proof of Theorem 11.** Define $g_j(t) = n_k^{-1} \sum z_{ik}(x_{ij} - t\hat{\mu}_j)\|x_i - t\hat{\mu}\|_{\hat{\Sigma}}^{-1}$, where $0 \leq t \leq 1$ and $j = 1, \ldots, p$. Then the $j$th component of $n_k^{-1} \sum z_{ik} \times (x_i - \hat{\mu})/\|x_i - \hat{\mu}\|_{\hat{\Sigma}}$ is $g_j(1)$. Because the first-order derivative of $g_j(t)$, $\dot{g}_j(\cdot)$, is a continuous function of $t$, there must exist a $0 \leq \tilde{t} \leq 1$ such that $g_j(1) - g_j(0) = \dot{g}_j(\tilde{t})$, where

$$\begin{aligned}
\dot{g}_j(t) = &-\mu_j\Big\{n_k^{-1}\sum z_{ik}\|x_i - t\hat{\mu}\|_{\hat{\Sigma}}^{-1}\Big\} \\
&+ n_k^{-1}\sum z_{ik}(x_{ij} - t\hat{\mu}_j)\|x_i - t\hat{\mu}\|_{\hat{\Sigma}}^{-3}(x_i - t\hat{\mu})^\top \hat{\Sigma}^{-1}\hat{\mu}.
\end{aligned} \quad \text{(A.24)}$$

Let $\tilde{\mu} = \tilde{t}\hat{\mu} = (\tilde{\mu}_1, \ldots, \tilde{\mu}_p)^\top \in \mathbb{R}^p$. Then, by (A.24), we obtain

$$\begin{aligned}
&\left|n_k^{-1}\sum_{i=1}^n z_{ik}\left\{\frac{x_{ij} - \hat{\mu}_j}{\|x_i - \hat{\mu}\|_{\hat{\Sigma}}} - \frac{x_{ij}}{\|x_i\|_{\hat{\Sigma}}}\right\}\right| \\
&= |g_j(1) - g_j(0)| = |\dot{g}_j(\tilde{t})| \\
&= \left|n_k^{-1}\sum_{i=1}^n z_{ik}\left\{\frac{-e_j}{\|x_i - \tilde{\mu}\|_{\hat{\Sigma}}} + \frac{(x_{ij} - \tilde{\mu}_j)\hat{\Sigma}^{-1}(x_i - \tilde{\mu})}{\|x_i - \tilde{\mu}\|_{\hat{\Sigma}}^3}\right\}^\top \hat{\mu}\right|
\end{aligned} \quad \text{(A.25)}$$

for $1 \leq j \leq p$, where $e_j$ is defined in Section 2.3. After simple calculations, the right-hand side of (A.25) is bounded by

$$n_k^{-1}\sum_{i=1}^n \left\{\frac{|\hat{\mu}_j|}{\|x_i - \tilde{\mu}\|_{\hat{\Sigma}}} + \frac{|x_{ij} - \tilde{\mu}_j| \times |\hat{\mu}^\top \hat{\Sigma}^{-1}(x_i - \tilde{\mu})|}{\|x_i - \tilde{\mu}\|_{\hat{\Sigma}}^3}\right\}$$

$$\leq n_k^{-1}\sum_{i=1}^n \left\{\frac{\|\hat{\mu}\|}{\|\hat{\Sigma}\|^{-1/2}\|x_i - \tilde{\mu}\|} + \frac{\|x_i - \tilde{\mu}\| \times \|\hat{\mu}\|_{\hat{\Sigma}} \times \|x_i - \tilde{\mu}\|_{\hat{\Sigma}}}{\|\hat{\Sigma}\|^{-3/2}\|x_i - \tilde{\mu}\|^3}\right\}$$

$$\leq n_k^{-1}\sum_{i=1}^n \left\{\frac{\|\hat{\mu}\|}{\|\hat{\Sigma}\|^{-1/2}\|x_i - \tilde{\mu}\|} + \frac{\|\hat{\Sigma}^{-1}\| \times \|\hat{\mu}\|}{\|\hat{\Sigma}\|^{-3/2}\|x_i - \tilde{\mu}\|}\right\} \quad \text{(A.26)}$$

$$= \frac{n}{n_k} \times \|\hat{\mu}\| \times \|\hat{\Sigma}\|^{1/2}$$

$$\times (1 + \|\hat{\Sigma}^{-1}\| \times \|\hat{\Sigma}\|) \times \left\{ n^{-1} \sum_{i=1}^{n} \frac{1}{\|x_i - \tilde{\mu}\|} \right\}.$$

Obviously, $n_k/n = O_p(1)$. Because $\hat{\Sigma} \to_p \Sigma$, we have both $\|\hat{\Sigma}\|^{1/2} = O_p(1)$ and $\|\hat{\Sigma}^{-1/2}\| = O_p(1)$. These results together with $\|\hat{\mu}\| = O_p(n^{-1/2})$ imply that the right-hand side of (A.26) is $O_p(n^{-1/2})$ if we are able to show that $n^{-1} \sum \|x_i - \tilde{\mu}\|^{-1} = O_p(1)$. To this end, we compute

$$\left| n^{-1} \sum_{i=1}^{n} \|x_i - \tilde{\mu}\|^{-1} - n^{-1} \sum_{i=1}^{n} \|x_i\|^{-1} \right|$$

(A.27)
$$\leq n^{-1} \sum_{i=1}^{n} \|x_i\|^{-1} \times \left| \frac{\|x_i\|}{\|x_i - \tilde{\mu}\|} - 1 \right|$$

$$\leq \left( \max_{1 \leq i \leq n} \left| \frac{\|x_i\|}{\|x_i - \tilde{\mu}\|} - 1 \right| \right) \times \left( n^{-1} \sum_{i=1}^{n} \|x_i\|^{-1} \right).$$

By Theorem 10(ii), the first term on the right-hand side of (A.27) is $o_p(1)$. In addition, by Theorem 10(i) and the Law of Large Numbers, the second term on the right-hand side of (A.27) is $O_p(1)$. As a result, the right-hand side of (A.27) is $o_p(1)$. This together with the result of $n^{-1} \sum_{i=1}^{n} \|x_i\|^{-1} = O_p(1)$ means that the last term on the right-hand side of (A.26) is $O_p(1)$. Hence, the left-hand side of (A.25) is $O_p(n^{-1/2})$. Next, consider

(A.28)
$$\left| n_k^{-1} \sum_{i=1}^{n} z_{ik} \left\{ \frac{x_{ij}}{\|x_i\|_{\hat{\Sigma}}} - \frac{x_{ij}}{\|x_i\|} \right\} \right|$$

$$= \left| n_k^{-1} \sum_{i=1}^{n} z_{ik} \left\{ \frac{x_{ij}}{\|x_i\|_{\hat{\Sigma}} \times \|x_i\|} \times \frac{\|x_i\|^2 - \|x_i\|_{\hat{\Sigma}}^2}{\|x_i\| + \|x_i\|_{\hat{\Sigma}}} \right\} \right|$$

$$= \left| n_k^{-1} \sum_{i=1}^{n} z_{ik} \left\{ \frac{x_{ij}}{\|x_i\|_{\hat{\Sigma}} \times \|x_i\|} \times \frac{x_i^\top (I_p - \hat{\Sigma}^{-1}) x_i}{\|x_i\| + \|x_i\|_{\hat{\Sigma}}} \right\} \right|$$

$$\leq n_k^{-1} \sum_{i=1}^{n} \frac{\|x_i\|}{\|x_i\|_{\hat{\Sigma}} \times \|x_i\|} \times \frac{\|I_p - \hat{\Sigma}^{-1}\| \times \|x_i\|^2}{\|x_i\| + \|x_i\|_{\hat{\Sigma}}}$$

$$\leq \|I_p - \hat{\Sigma}^{-1}\| \times n_k^{-1} \sum_{i=1}^{n} \left| \frac{\|\hat{\Sigma}\|^{1/2} \times \|x_i\|}{\|x_i\| \times \|x_i\|} \times \frac{\|x_i\|^2}{\|x_i\| + \|\hat{\Sigma}\|^{-1/2} \|x_i\|} \right|$$

$$= \|I_p - \hat{\Sigma}^{-1}\| \times (n_k^{-1} n) \times \|\hat{\Sigma}\|^{1/2} \{1 + \|\hat{\Sigma}\|^{-1/2}\}^{-1}.$$

By Theorem 4.2 of Tyler (1987), we have $\|I_p - \hat{\Sigma}^{-1}\| = O_p(n^{-1/2})$. This implies that $\|\hat{\Sigma}\|^{-1/2} = O_p(1)$. Furthermore, $n_k^{-1} n = O_p(1)$. Accordingly, the right-hand side of (A.28) is of order $O_p(n^{-1/2})$. This indicates that the left-hand side of (A.28) is also of order $O_p(n^{-1/2})$. This in conjuction with (A.25) implies that

$$(A.29) \qquad n_k^{-1} \sum_{i=1}^n z_{ik} \left\{ \frac{x_i - \hat{\mu}}{\|x_i - \hat{\mu}\|_{\hat{\Sigma}}} - \frac{x_i}{\|x_i\|} \right\} = O_p(n^{-1/2}).$$

Because the $L_2$-norm of $x_i/\|x_i\|$ is 1, we then apply the central limit theorem to have $n_k^{-1} \sum z_{ik} x_i/\|x_i\| = E(\overrightarrow{X}|Y=k) + O_p(n^{-1/2})$. This result together with (A.29) shows that $n_k^{-1} \sum z_{ik}(x_i - \hat{\mu})/\|x_i - \hat{\mu}\|_{\hat{\Sigma}} = E(\overrightarrow{X}|Y=k) + O_p(n^{-1/2})$. Finally, using the fact that $\hat{\Sigma} - I_p = O_p(n^{-1/2})$, we have

$$\bar{x}_k = \hat{\Sigma}^{-1/2} \frac{1}{n_k} \sum z_{ik} \left( \frac{x_i - \hat{\mu}}{\|x_i - \hat{\mu}\|_{\hat{\Sigma}}} \right) = E(\overrightarrow{X}|Y=k) + O_p(n^{-1/2}).$$

This completes the proof. $\square$

**Acknowledgment.** We are very grateful to the Editor, the Associate Editor and the referees for their constructive comments and insightful suggestions that enrich the manuscript substantially.

R. Luo  
School of Finance  
Southwestern University of Finance  
and Economics  
Chengdu, Si Chuan Province  
P.R. China, 610074  
E-mail: ronghua@swufe.edu.cn

H. Wang  
Guanghua School of Management  
Peking University  
Beijing  
P.R. China, 100871  
E-mail: hansheng@gsm.pku.edu.cn

C.-L. Tsai  
Graduate School of Management  
University of California  
Davis, 95616-8609  
USA  
E-mail: cltsai@ucdavis.edu